
\def\LaTeX{%
  \let\Begin\begin
  \let\End\end
  \let\salta\relax
  \let\finqui\relax
  \let\futuro\relax}

\def\UK{\def\our{our}\let\sz s}
\def\USA{\def\our{or}\let\sz z}

\UK



\LaTeX

\USA


\salta

\documentclass[twoside,12pt]{article}
\setlength{\textheight}{24cm}
\setlength{\textwidth}{16cm}
\setlength{\oddsidemargin}{2mm}
\setlength{\evensidemargin}{2mm}
\setlength{\topmargin}{-15mm}
\parskip2mm


\usepackage[usenames,dvipsnames]{color}
\usepackage{amsmath}
\usepackage{amsthm}
\usepackage{amssymb}
\usepackage[mathcal]{euscript}
%
%
\usepackage{cite}
%
%


\definecolor{viola}{rgb}{0.3,0,0.7}
\definecolor{ciclamino}{rgb}{0.5,0,0.5}
\definecolor{rosso}{rgb}{0.8,0,0}
\definecolor{someblue}{rgb}{0.54, 0.81, 0.84}

\def\pier #1{{\color{rosso}#1}}
\def\pcol #1{{\color{rosso}#1}}
\def\lastrev #1{{\color{rosso}#1}}
\def\andrea #1{{\color{rosso}#1}}
\def\dafare #1{{\color{rosso}#1}}

\def\pier #1{#1}
\def\pcol #1{#1}
\def\lastrev #1{{#1}}
\def\andrea #1{#1}
\def\dafare  #1{#1}


\bibliographystyle{plain}


%

\finqui

\def\Beq{\Begin{equation}}
\def\Eeq{\End{equation}}
\def\Bsist{\Begin{eqnarray}}
\def\Esist{\End{eqnarray}}

\def\Bthm{\Begin{theorem}}
\def\Ethm{\End{theorem}}
\def\Blem{\Begin{lemma}}
\def\Elem{\End{lemma}}
\def\Bprop{\Begin{proposition}}
\def\Eprop{\End{proposition}}

\def\Brem{\Begin{remark}\rm}
\def\Erem{\End{remark}}

\def\Bdim{\Begin{proof}}
\def\Edim{\End{proof}}
\def\Bcenter{\Begin{center}}
\def\Ecenter{\End{center}}
\let\non\nonumber




\def\step #1 \par{\medskip\noindent{\bf #1.}\quad}


\def\Lip{Lip\-schitz}
\def\Holder{H\"older}
\def\Frechet{Fr\'echet}
\def\Poincare{Poincar\'e}
\def\aand{\quad\hbox{and}\quad}

\def\wk{well-known}
\def\socal{so-called}
\def\lhs{left-hand side}
\def\rhs{right-hand side}
\def\sfw{straightforward}

\def\CH{Cahn--Hilliard}


\def\generaliz{generali\sz}


\def\multibold #1{\def\arg{#1}%
  \ifx\arg\pto \let\next\relax
  \else
  \def\next{\expandafter
    \def\csname #1#1#1\endcsname{{\bf #1}}%
    \multibold}%
  \fi \next}

\def\pto{.}

\def\multical #1{\def\arg{#1}%
  \ifx\arg\pto \let\next\relax
  \else
  \def\next{\expandafter
    \def\csname cal#1\endcsname{{\cal #1}}%
    \multical}%
  \fi \next}


\def\multimathop #1 {\def\arg{#1}%
  \ifx\arg\pto \let\next\relax
  \else
  \def\next{\expandafter
    \def\csname #1\endcsname{\mathop{\rm #1}\nolimits}%
    \multimathop}%
  \fi \next}

\multibold
qwertyuiopasdfghjklzxcvbnmQWERTYUIOPASDFGHJKLZXCVBNM.

\multical
QWERTYUIOPASDFGHJKLZXCVBNM.

\multimathop
ad dist div dom meas sign supp .


\def\accorpa #1#2{\eqref{#1}--\eqref{#2}}
\def\Accorpa #1#2 #3 {\gdef #1{\eqref{#2}--\eqref{#3}}%
  \wlog{}\wlog{\string #1 -> #2 - #3}\wlog{}}


\def\neto{\mathrel{{\scriptscriptstyle\nearrow}}}
\def\seto{\mathrel{{\scriptscriptstyle\searrow}}}

\def\graffe #1{\mathopen\{#1\mathclose\}}

\def\<#1>{\mathopen\langle #1\mathclose\rangle}
\def\norma #1{\mathopen \| #1\mathclose \|}

\def\iot {\int_0^t}
\def\ioT {\int_0^T}
\def\intQt{\int_{Q_t}}
\def\intQ{\int_Q}
\def\iO{\int_\Omega}
\def\iG{\int_\Gamma}
\def\intS{\int_\Sigma}
\def\intSt{\int_{\Sigma_t}}

\def\dt{\partial_t}
\def\dn{\partial_n}
\def\dtt{\partial_{tt}}

\def\cpto{\,\cdot\,}

\def\checkmmode #1{\relax\ifmmode\hbox{#1}\else{#1}\fi}

\def\aeQ{\checkmmode{a.e.\ in~$Q$}}

\def\aeS{\checkmmode{a.e.\ on~$\Sigma$}}
\def\aet{\checkmmode{a.e.\ in~$(0,T)$}}

\def\aat{\checkmmode{for a.a.~$t\in(0,T)$}}


\def\erre{{\mathbb{R}}}




\def\genspazio #1#2#3#4#5{#1^{#2}(#5,#4;#3)}
\def\spazio #1#2#3{\genspazio {#1}{#2}{#3}T0}

\def\L {\spazio L}
\def\H {\spazio H}
\def\W {\spazio W}

\def\C #1#2{C^{#1}([0,T];#2)}


\def\Lx #1{L^{#1}(\Omega)}
\def\Hx #1{H^{#1}(\Omega)}

\def\LxG #1{L^{#1}(\Gamma)}
\def\HxG #1{H^{#1}(\Gamma)}

\def\LS #1{L^{#1}(\Sigma)}

\def\Ldue{\Lx 2}

\def\Huno{\Hx 1}
\def\Hdue{\Hx 2}

\def\HunoG{\HxG 1}
\def\HdueG{\HxG 2}

\def\LdueG{\LxG 2}



\let\theta\vartheta

\let\phi\varphi
\let\lam\lambda

\let\TeXchi\chi                         
\newbox\chibox
\setbox0 \hbox{\mathsurround0pt $\TeXchi$}
\setbox\chibox \hbox{\raise\dp0 \box 0 }
\def\chi{\copy\chibox}



\def\suG{{\vrule height 5pt depth 4pt\,}_\Gamma}

\def\fG{f_\Gamma}
\def\yG{y_\Gamma}
\def\uG{u_\Gamma}
\def\vG{v_\Gamma}
\def\hG{h_\Gamma}
\def\xiG{\xi_\Gamma}
\def\qG{q_\Gamma}
\def\uGmin{u_{\Gamma,{\rm min}}}
\def\uGmax{u_{\Gamma,{\rm max}}}

\def\yz{y_0}

\def\Mz{M_0}
\def\ustar{u_*}
\def\vstar{v_*}

\def\bQ{b_Q}
\def\bS{b_\Sigma}
\def\bO{b_\Omega}
\def\bG{b_\Gamma}
\def\bz{b_0}

\def\rmin{r_-}
\def\rmax{r_+}

\def\zQ{z_Q}
\def\zS{z_\Sigma}
\def\zO{z_\Omega}
\def\zG{z_\Gamma}
\def\phQ{\phi_Q}
\def\phS{\phi_\Sigma}
\def\phO{\phi_\Omega}
\def\phG{\phi_\Gamma}

\def\Uad{\calU_{\ad}}
\def\uopt{\overline u_\Gamma}
\def\yopt{\overline y}
\def\yGopt{\overline y_\Gamma}
\def\wopt{\overline w}

\def\redJ{\widetilde\calJ}

\def\unG{u_{\Gamma,n}}
\def\yn{y_n}
\def\ynG{y_{\Gamma,n}}
\def\wn{w_n}

\def\yh{y^h}
\def\yhG{y^h_\Gamma}
\def\wh{w^h}
\def\qh{\theta^h}
\def\qhG{\theta^h_\Gamma}
\def\zh{z^h}

\def\uO{u^\Omega}
\def\vO{v^\Omega}

\def\VO{\calV_\Omega}

\def\GO{\calG_\Omega}

\def\VG{V_\Gamma}
\def\HG{H_\Gamma}
\def\Vp{V^*}

\def\normaV #1{\norma{#1}_V}
\def\normaH #1{\norma{#1}_H}
\def\normaVG #1{\norma{#1}_{\VG}}

\def\normaVp #1{\norma{#1}_*}

\def\nablaG{\nabla_{\!\Gamma}}
\def\DeltaG{\Delta_\Gamma}

\let\hat\widehat

\def\Pi{\hat\pi}

\def\lamG{\lam_\Gamma}

\def\mz{m_0}


\def\cd{c_\delta}

\Begin{document}



\title{Boundary control problem and\\[0.3cm] 
  optimality conditions\\[0.3cm]
  for the \CH~equation\\[0.3cm]
  with dynamic boundary conditions}
\author{}
\date{}
\maketitle

\Bcenter
\vskip-1cm
{\large\sc Pierluigi Colli$^{(1)}$}\\
{\normalsize e-mail: {\tt pierluigi.colli@unipv.it}}\\[.25cm]
{\large\sc Andrea Signori$^{(2)}$}\\
{\normalsize e-mail: {\tt andrea.signori02@universitadipavia.it}}\\[.25cm]
$^{(1)}$
{\small Dipartimento di Matematica ``F. Casorati'', Universit\`a di Pavia}\\
{\small via Ferrata 5, 27100 Pavia, Italy}\\[.2cm]
$^{(2)}$
{\small Dipartimento di Matematica e Applicazioni, Universit\`a di Milano--Bicocca}\\
{\small via Cozzi 55, 20125 Milano, Italy}

\Ecenter
\Begin{abstract}\noindent
This paper is concerned with a boundary control problem for the 
Cahn--Hilliard equation coupled with dynamic boundary conditions.
In order to handle the control problem,
we restrict our analysis to the case of regular
potentials defined on the whole real line,
assuming the boundary potential to be dominant.
The existence of optimal control, the Fr\'echet differentiability
of the control-to-state \andrea{operator} between appropriate Banach spaces, 
and the first-order necessary conditions for optimality 
are addressed.
In particular, the \andrea{necessary} condition for optimality is characterized by
a variational inequality involving the adjoint variables.

\vskip3mm
\noindent {\bf Key words:}
\CH~equation, dynamic boundary conditions, phase separation,
double-well potentials, optimal control, optimality conditions,
adjoint problem.
\vskip3mm
\noindent {\bf AMS (MOS) Subject Classification:} {35K61, 49J20, 49K20, 49J50.}
\End{abstract}

\pagestyle{myheadings}
\newcommand\testopari{\sc Colli -- Signori}
\newcommand\testodispari{\sc Boundary control problem for the \CH~equation}
\markboth{\testopari}{\testodispari}

\salta
\finqui

\newpage
\section{Introduction}
\label{Intro}
\setcounter{equation}{0}
The \CH~equation plays a fundamental role in material science
(see, e.g., the review paper \cite{Mir_CH} and the vast literature therein). 
Such an equation was historically proposed for the study of phase \andrea{segregation} in 
cooling binary alloys (see \cite{CH}). On the other hand, 
from then onward, it has been shown how versatile this
equation can be for several applications in very different fields such as engineering, 
biology, tumor growth,
image inpainting, population dynamics, bacterial films, and many others.
The huge efforts by the mathematical community have made the classical \CH~equation
well-understood from a mathematical point of view,
at least as far as the existence, uniqueness and
regularity of solutions are concerned. Here, we address a boundary optimal control problem 
for the \CH~equation coupled with some non-standard boundary conditions, the 
\socal~dynamic ones.

For a fixed finite final time $T > 0$, the \CH~equation reads as follows
\begin{align}
  \label{cahnh_uno}
  & \dt y - \Delta w = 0   \quad \hbox{in $Q:=\Omega\times (0,T) $},
  \\ 
  & w = -\Delta y + f'(y) 
  \quad \hbox{in $Q$},
  \label{cahnh_due}
\end{align}
\Accorpa\cahn cahnh_uno cahnh_due
where $\Omega$ represents the space domain in which the evolution takes place, 
and the occurring variables $y$~and $w$ stand for the order parameter 
and the corresponding chemical potential, respectively. 
Moreover, $f'$ denotes the derivative of a nonlinearity that possesses a double-well behavior. 
For this latter, \pier{the \lastrev{prototype} is the regular double-well potential 
$f_{reg}$, defined by 
\begin{align}
   f_{reg}(r) = \frac14(r^2-1)^2, \quad  \hbox{whence} \quad  
   f_{reg}'(r) = r^3- r,   \quad \hbox{for } r \in \erre.
  \label{regpot}
\end{align}
Besides,} we endow the above system with an initial condition 
of the form $y(0)=\yz$, and suitable boundary conditions. 
As \pier{for boundary conditions,} the widespread types in literature are the 
no-flux conditions for both the variables $y$ and~$w$. It is worth noting that,
from a phenomenological point of view, the \pier{no-flux condition for $w$}
is quite natural since it ensures the mass conservation
during the evolution process\pier{: this can} be easily 
checked by testing the equation 
\eqref{cahnh_uno} by $1$ and integrating \pier{by parts} over $\Omega$.
In fact, denoting by $(v)^{\Omega}$ the mean value of the function $v: \Omega \to \erre$,
we realize that
\begin{align}
  & (\dt y(t))^\Omega = 0
  \quad \aat,
  \aand
  (y(t))^\Omega = \mz
  \quad \hbox{for every $t\in[0,T]$,}
  \quad
  \non
  \\
  &\quad  \hbox{where $\mz:=(\yz)^\Omega$ is the mean value of $\yz$.}
  \label{conserved}
\end{align}
In this contribution, we also deal with the no-flux condition for
the chemical potential, whereas a dynamic boundary condition for the order parameter 
is \andrea{prescribed}. 
These boundary conditions are quite new and were recently proposed in order to 
take into account the dynamics between the walls. 
In this regard, let us address to \cite{CGS}, where both the 
viscous and the non-viscous \CH~equations, combined with these kinds of boundary 
conditions, have been investigated by assuming the boundary potential to be dominant
on the bulk one. Furthermore, we \pier{have to} mention\pier{\cite{CFP,CF2,GiMiSchi,PRZ,RZ, MZ, CGS_Lim, CGS_Rec,WZ,GK, LW}},
where other problems related to the \CH~equation combined with \pier{dynamic} boundary 
conditions have been analyzed\pier{,
and\pier{\cite{CaCo,CF, CFS,Is, CGS_nonpf, MRSS_pf, CGNS_AC}} for the coupling of dynamic boundary conditions with} different phase field 
models such as the Allen-Cahn or the Penrose-Fife model.
So, \pier{according to \cite{CGS} we supply the above system 
\eqref{cahnh_uno}--\eqref{cahnh_due} with}
\begin{align}
  & \dn w =0 
  \quad \hbox{on $\Sigma:=\Gamma\times (0,T)$},
  \label{nfc}
  \\ & 
  \dn y + \dt\yG - \DeltaG\yG + f'_\Gamma(\yG) = \uG
  \quad \hbox{on $\Sigma,$}
  \label{dbc}
\end{align}
where $\Gamma$ \pier{is} the boundary of $\Omega$,
$\yG$ denotes the trace of $y$, 
$\DeltaG$ stands for the Laplace--Beltrami operator on the boundary, and
$\dn$ \pier{represents} the outward normal derivative. Furthermore, the term
$f'_\Gamma$ is a nonlinearity similar to the previous $f'$,
but operating on the values on the boundary instead of on the bulk, whereas 
$\uG$ is the \socal~control variable which can be interpreted as a boundary source term. 

Summing up, the system we \andrea{are going} to deal with reads
\begin{align}
  & \dt y - \Delta w = 0
  \quad \hbox{in $\, Q$},
  \label{Iprima}
  \\
  & w = - \Delta y + f'(y)
  \quad \hbox{in $\,Q$},
  \label{Iseconda}
  \\
  & \dn w = 0
  \quad \hbox{on $\, \Sigma$},
  \label{Ibc}
  \\
  & \yG = y\suG
  \aand
  \dt\yG + \dn y - \DeltaG\yG + f'_\Gamma(\yG) = \uG
  \quad \hbox{on $\, \Sigma$},
  \label{Iterza}
  \\
  & y(0) = \yz
  \quad \hbox{in $\, \Omega$} .
  \label{Icauchy}
\end{align}
\Accorpa\Ipbl Iprima Icauchy
Once that the state system \Ipbl~has been described, we can address
the corresponding control problem. 
Among several possibilities, we consider the following tracking-type cost functional
\begin{align}
  \calJ(y,\yG,\uG)
  &:= \frac\bQ 2 \, \norma{y-\zQ}_{L^2(Q)}^2
  + \frac\bS 2 \, \norma{\yG-\zS}_{L^2(\Sigma)}^2
  + \frac\bO 2 \, \norma{y(T)-\zO}_{L^2(\Omega)}^2
  \non
  \\
  & \quad\quad
  + \frac\bG 2 \, \norma{\yG(T)-\zG}_{L^2(\Gamma)}^2
  + \frac\bz 2 \, \norma\uG_{L^2(\Sigma)}^2,
  \label{defcost}
\end{align}
where the symbols $\bQ,\bS,\bO,\bG,\bz$ and $z_Q, z_\Sigma,z_\Omega,z_\Gamma$ denote
nonnegative constants and some target functions, respectively.
Moreover, we require the control variable $\uG$ to belong to the 
non-empty control-box 
$\Uad$ which is defined by
\begin{align}
  \Uad & := 
  \bigl\{ \uG\in\H1 {L^2(\Gamma)}\cap\LS\infty:
  \non
  \\
  & \qquad\quad
  \uGmin\leq\uG\leq\uGmax\ \ \aeS,\ \, \norma{\dt\uG}_{L^2(\Sigma)}\leq\Mz
  \bigr\},
  \label{defUad}
\end{align}
for suitable functions $\uGmin,\uGmax\in\LS\infty$, and for a positive constant $\Mz$.
Note that, owing to the weak lower semicontinuity of norms, $\Uad$ is a 
closed convex subset of $L^2(\Sigma)$.
Therefore, our minimization problem consists in seeking
an admissible control variable $\uG$ such that, along with its corresponding solution
to system \Ipbl, minimizes the cost functional \eqref{defcost}.

Concerning the interpretation of the optimal control problem, 
let us point out that, since the target functions 
$z_Q, z_\Sigma,z_\Omega,z_\Gamma$ \pier{provide some
particular configuration}, we are looking for 
an admissible control variable $\uG$ which \pier{forces
its corresponding solution to \Ipbl\ to be} as close as 
possible to the prescribed configuration.
Conversely, the last term of \eqref{defcost} penalizes the large values
of the \pier{$L^2$-norm} of the control \andrea{so that} 
it can be seen as the cost we have to pay in order to follow
that strategy.

\lastrev{As for previous contributions on optimal control problems 
\pcol{for \CH\ systems possibly involving
dynamic boundary conditions, let us mention the papers}~\cite{CFGS1,CFGS2,
CGM,GS,CGS_vel, CGS_nonst,CGS_OPT,CGSpure,CS_opt,CGSGuid,FY,
Hint_uno,Hint_due,RS_OPT,ZL,ZL_due}.
In particular, we focus \pcol{our} attention on~\cite{CGS_OPT}, 
where the optimal control problem for the viscous \CH~equation 
endowed with dynamic boundary conditions is investigated by
exploiting the well-posedness of the state system 
discussed in \cite{CGS}.}
Moreover, we \lastrev{also point out} \cite{CGSpure}, where 
the optimal control problem \pier{has been extended} to the non-viscous case. 
In fact, by employing suitable asymptotic arguments \pier{and letting the viscosity parameter \pier{tend to} zero, \pier{in \cite{CGSpure} it is shown}} how 
the optimal control \pier{results for the viscous case 
allow to recover other} results for the pure setting.
It is worth underlining that the optimal control
problem is exactly the one we are going to address here\pier{,
but in this contribution we follow a direct approach and are able to obtain better results.}

\andrea{\pier{Indeed, it occurs that in the limit procedure of \cite{CGSpure} 
some information on the limiting terms turns out to be 
lost and especially the results concerning the first-order conditions for
optimality and the adjoint system are somehow unsatisfactory since they
hold in a very weak sense. Moreover,
the adjoint system at the limit has not an explicit structure
and the uniqueness for its solution is not at all clear.}
Namely, the related \pier{existence  result 
states the existence of proper elements in dual spaces that
satisfy some properties and are the (weak star) limits of
some terms or groups of terms of the adjoint system for the viscous case 
(see \cite[Thm.~2.7, p.~318]{CGSpure} for \pier{a} precise statement).
On the other hand, it may appear that the optimality condition there obtained
is very similar to the one we will point out here since
they formally consist in the same variational inequality}
\lastrev{(see \cite[eq.~(2.54)]{CGSpure} and compare with \eqref{cnoptadj}).}
\pier{However, the results are substantially different and the 
difference is hidden in the two adjoint systems}.
Lastly, let us point out that the cost functional
of \cite{CGSpure} is less general than \pier{ours} since,
as a consequence of the results of \cite{CGS_OPT}, the constant $b_\Omega$
and $b_\Gamma$ are taken identically zero.}

\andrea{\pier{In the present contribution, provided} we restrict the analysis on everywhere defined potentials like \pier{\eqref{regpot}, we show} that also for the non-viscous case
the optimality condition can be completely characterized.
As a matter of fact, the existence, uniqueness and also \pier{further} regularity
for the adjoint system will be proved (\pier{cf.}~Theorem~\ref{adjexist}).
Moreover, since from a technical viewpoint the strategies are very different\pier{,}
we have to perform the proofs ex novo without relying on the results
proved in \cite{CGS_OPT}.}

\andrea{After showing the existence of optimal controls,
\pier{we characterize} the} first-order necessary 
conditions that every optimal control has to satisfy through a variational inequality.
In this direction, a key point will be showing the \Frechet\ differentiability of the
control-to-state \andrea{operator}. Then, as usual for optimal control problems 
(see, e.g., \cite{Trol, Lions_OPT}), 
in order to simplify
the obtained optimality conditions, a new system, called adjoint, 
has to be introduced and solved in order to reformulate the necessary condition
in a more convenient way.
The adjoint system turns out to be a backward-in-time boundary value problem 
of the following form
\begin{align}
  \non
  & q = -\Delta p  \quad \hbox{in $Q$},
  \\ &\non
  - \dt p - \Delta q + \lam q = \phi_Q
  \quad \hbox{in $Q$},
  \\ & \non
  \dn p = 0 \quad \hbox{on $\Sigma$},
  \\ &\non
  - \dt\qG + \dn q - \DeltaG\qG + \lamG \qG = \phi_\Sigma
  \quad \hbox{on $\Sigma$},
\end{align}
where \pier{$q$ and $p$ are the adjoint variables,} 
$\qG$ stands for the trace of $q$, \pier{and} the functions $\lam$, $\lamG$, $\phi_Q$ 
and $\phi_\Sigma$ are somehow related to 
$z_Q, z_\Sigma,z_\Omega,z_\Gamma$ and to the constants $\bQ,\bS,\bO,\bG,\bz$ 
appearing in \eqref{defcost}, as well as to the optimal state $(\yopt,\yGopt)$, 
which is the state associated to the optimal control $\uopt$. 
Furthermore, the above system will be coupled with suitable final conditions.

\pier{The plan of the paper} is as follows.
In Section \ref{STATEMENT} we \pier{specify} the mathematical setting 
and recollect the results we have established.
From the third section on, we begin with the corresponding proofs. 
Section \ref{OPTIMUM} is devoted \pier{to the existence of optimal controls}.
Furthermore, Section \ref{FRECHET} is the place in which 
the main novelties \andrea{appear\pier{: there}, we} discuss the properties of the control-to-state
\andrea{operator} $\cal S$ proving its \Lip~continuity and 
the \Frechet~differentiability in suitable Banach spaces. 
Finally, the \andrea{well-posedness of the}
adjoint system and the first-order necessary conditions 
for optimality \pier{are discussed} \andrea{in Section \ref{OPTIMALITY}.}

\section{Statement of the problem and results}
\label{STATEMENT}
\setcounter{equation}{0}
In this section, we set \andrea{the} notation \pier{and present} 
in detail \andrea{the established} results.
\pier{We start by pointing} out that $\Omega$~represents the 
body where the evolution takes place 
and we assume $\Omega\subset\erre^3$ to be open, connected, bounded and smooth,
with Lebesgue measure denoted by $|\Omega|$. 
Moreover, let us fix once for all that the 
symbols $\Gamma$, $\dn$, $\nablaG$ and $\DeltaG$ stand for the boundary of~$\Omega$, 
the outward normal derivative, the surface gradient, and the Laplace--Beltrami operator, 
respectively.
Given a finite final time~$T>0$,
we set for convenience
\begin{align}
  & Q_t := \Omega \times (0,t),
  \quad
  \Sigma_t := \Gamma \times (0,t)
  \quad \hbox{for every $t\in(0,T]$},
  \label{defQtSt}
  \\
  & Q := Q_T \,,
  \aand
  \Sigma := \Sigma_T \,.
  \label{defQS}
\end{align}
\Accorpa\defQeS defQtSt defQS
Before diving into the mathematical setting, let us
emphasize a typical issue of control problems. Although 
some of the results we need hold under rather weak conditions,
we will require quite strong hypotheses for the involved
potentials and for the initial data
in order to handle the corresponding control problem.
As a consequence, the following requirements surely comply with the framework of \cite{CGS}.
 
On the potentials \andrea{$f$ and $\fG$} we make the following structural assumptions
\begin{align}
  & \hbox{$f,\,\fG: \erre\to[0,+\infty)$
   are {$C^4$} functions.}
  \label{hppot}
  \\
  & f'(0) = \fG'(0) = 0,
  \, \hbox{and }
  \hbox{$f''$ and $\fG''$ are bounded from below.}
  \label{hpfseconda}
  \\
  & |f'(r)| \leq \eta \,|\fG'(r)| + C
  \quad \hbox{for some $\eta,\, C>0$ \, and every $r\in\erre$.}
  \label{hpcompatib}
  \\
  & \lim\limits_{r\seto\ -\infty} f'(r)
  = \lim\limits_{r\seto -\infty} \fG'(r) 
  = -\infty ,
  \,\, \hbox{and }
  \lim\limits_{r\neto +\infty} f'(r)
  = \lim\limits_{r\neto +\infty} \fG'(r) 
  = +\infty \, .
  \label{fmaxmon}
\end{align}
\Accorpa\HPstruttura hppot fmaxmon
\Brem
The above conditions imply the possibility of splitting $f'$ as $f'=\beta+\pi$,
where $\beta$ is a monotone function, which diverges as its argument goes to $-\infty$
or to $+\infty$, while $\pi$ is a regular perturbation with bounded derivative. 
\andrea{Likewise}, it goes for the boundary contribution~$\fG'$ that can be 
\andrea{possibly} written as 
$\fG'=\beta_\Gamma+\pi_\Gamma$, for suitable functions \pier{satisfying the same properties as $\beta$ and $\pi$}.
\Erem

It is worth emphasizing that in our treatment, owing to \HPstruttura, 
the case of \eqref{regpot} is allowed, while other significant cases \pier{like, e.g.,
the logarithmic potential 
\begin{align}
   f_{log}(r) &= (1+r)\ln (1+r)+(1-r)\ln (1-r) - k r^2 , 
  \quad r \in (-1,1),
  \label{logpot}
\end{align}
(with $k>1$ to ensure non-convexity)} \pier{are} not.
\pier{On the other hand, the above setting~\HPstruttura} perfectly fits the framework 
of \cite{CGS} since the assumption \eqref{hpcompatib} postulates
the domination of the boundary potential on the bulk one.
For the converse case, namely the one in which the bulk potential is the 
leading one between the two, \pier{we refer to} the contributions \cite{GiMiSchi,GMS_long}.
Now, let us introduce some functional spaces that will be useful later on
by defining
\begin{align}
   V &:= \Huno, \quad
  H := \Ldue, \quad
 \VG := \HunoG ,
  \quad
   \HG := \LdueG ,
  \qquad
  \label{defVH}
  \\
   \calV &:= \graffe{(v,\vG) \in V\times\VG:\ \vG=v\suG}, \aand
  \calG := \Vp \times \HG,
  \label{defcalVH}
\end{align}
\Accorpa\Defspazi defVH defcalVH
and we endow them with their natural norms to get some Banach spaces. 
Besides, for an arbitrary Banach space $X$, we agree to use
$\norma{\cdot}_{X}$ to denote its norm,
the standard symbol $X^*$ for its topological dual, and ${}_{X^*}\<\cdot, \cdot >_X$ 
for the corresponding duality product between $X^*$ and $X$. Meanwhile, we will
use $\norma\cpto_p$ for the usual norm in $L^p$ spaces. 
In the following, we understood that $H$ is embedded in $\Vp$
in the usual way, i.e. $V \subset H \cong H^* \subset \Vp$. This constitutes a
Hilbert triplet, namely we have the following identification
\Beq
   \label{ternahilbert}
   \<u,v> \ = (u,v) \quad{} \hbox{for every $u\in H$ and $v\in V$,}
\Eeq
where $(\,\cdot\,,\,\cdot)$ denotes the inner product in $H$. 

In addition, whenever $u\in\Vp$ and $\underline u\in\L1\Vp$,
we define their \generaliz ed mean values 
$\uO\in\erre$ and $\underline u^\Omega\in L^1(0,T)$ by
\Beq
  \uO := \frac 1 {|\Omega|} \, \< u , 1 >,
  \aand
  \underline u^\Omega(t) := \bigl( \underline u(t) \bigr)^\Omega
  \quad \aat,
  \label{media}
\Eeq
where \eqref{media} reduces to the usual mean values when it is applied to elements
of~$H$ or $\L1H$.

Next, since in the last two sections we are going to use
test functions with zero mean value, it is convenient to set
\Beq
  \GO := \graffe{(v,\vG)  \in \calG :\ \vO=0},
  \aand
  \VO := \GO \cap \calV,
  \label{defGOVO}
\Eeq
and endow them with their natural topologies as subspaces of $\calG$ and~$\calV$, respectively.
Moreover, we define
\Beq
  \dom\calN := \graffe{\vstar\in\Vp: \ \vstar^\Omega = 0},
  \aand
  \calN : \dom\calN \to \graffe{v \in V : \ \vO = 0},
  \label{predefN}
\Eeq
as the map which assigns to every $\vstar\in\dom\calN$ the element 
${\cal N} \vstar$ which satisfies
\Beq
  {\calN\vstar \in V, \quad
  (\calN\vstar)^\Omega = 0 ,
  \aand
  \iO \nabla\calN\vstar \cdot \nabla z = \< \vstar , z >
  \quad \hbox{for every $z\in V$}}.
  \label{defN}
\Eeq
Hence, $\calN\vstar$
represents the solution $v$ to the \generaliz ed Neumann problem for $-\Delta$
with datum~$\vstar$ that in addition has to satisfy the zero mean value condition.
In fact, if $\vstar\in H$, the \andrea{above} conditions mean that
$-\Delta\calN\vstar = \vstar$ in $\Omega$ and $\dn(\calN\vstar )= 0$ on $\Gamma$.
As far as $\Omega$ is bounded, smooth and connected,
it follows that \eqref{defN} yields a well-defined isomorphism
which also satisfies
\Bsist
  && \calN\vstar\in\Hx{s+2},
  \quad
  \norma{\calN\vstar}_{\Hx{s+2}}
  \leq C_s \norma\vstar_{\Hx s},
  \non
  \\
  && \quad \hbox{if $s\geq0$}
  \aand
  \vstar \in \Hx s \cap \dom\calN ,
  \label{regN}
\Esist
with a constant $C_s$ that depends only on $\Omega$ and~$s$.
Moreover, we have the following properties
\Beq
  \< \ustar , \calN \vstar >
  = \< \vstar , \calN \ustar >
  = \iO (\nabla\calN\ustar) \cdot (\nabla\calN\vstar)
  \quad \hbox{for $\ustar,\vstar\in\dom\calN$},
  \label{simmN}
\Eeq
whence also
\Beq
  2 \< \dt\vstar(t) , \calN\vstar(t) >
  = \frac d{dt} \iO |\nabla\calN\vstar(t)|^2
  = \frac d{dt} \, \normaVp{\vstar(t)}^2
  \quad \aat,
  \label{dtcalN}
\Eeq
for every $\vstar\in\H1\Vp$ satisfying $(\vstar)^\Omega=0$ \aet,
where we have set 
$\norma{\cdot}_* := \normaH{\nabla \calN (\cdot)}$\pier{, which 
turns out to be a norm in $\Vp$ equivalent to the usual dual norm.}

As the initial data are concerned, we require that
\begin{align}
  & \yz \in \Hx 2, \quad
  \yz\suG \in \HxG 2,
  \aand \Delta\yz \in V,
  \label{hpyz}
\end{align}
where the last condition has been already assumed in~\cite{CGS} to
ensure good regularity results for the non-viscous system (see \cite[Eq.~(2.40), p.~978]{CGS}).
Even though we could write the equations and the boundary conditions in their strong forms,
we however prefer to use the corresponding variational formulations. Hence, the problem we want 
to deal with consists of looking for a triplet $(y,\yG,w)$ that satisfies the regularity 
\begin{align}
 y  & \in \W{1,\infty}\Vp \cap \H1V \cap \L\infty\Hdue,
  \label{regy}
  \\
   \yG &\in \W{1,\infty}\HG \cap \H1\VG \cap \L\infty\HdueG,
  \label{regyG}
  \\
   \yG(t) &= y(t)\suG
  \quad \aat,
  \label{tracciay}
  \\
   w &\in \L\infty V \cap \L2 {\Hx3},
  \label{regw}
\end{align}
\Accorpa\Regsoluz regy regw
as well as, for almost every $t\in (0,T)$, the variational equalities
\begin{align}
  & \< \dt y(t) , v  >
  + \iO \nabla w(t) \cdot \nabla v = 0
  \quad \hbox{for every $v\in V$},
  \label{prima}
  \\
   & \iO w(t) \, v
  = \iG \dt\yG(t) \, \vG
  + \iO \nabla y(t) \cdot \nabla v
  + \iG \nablaG\yG(t) \cdot \nablaG\vG
  \non
  \\
  & \quad {}
  + \iO f'(y(t)) \, v
  + \iG \bigl( \fG'(\yG(t)) - \uG(t) \bigr) \, \vG
  \quad \hbox{for every $(v,\vG) \in\calV$}
  \label{seconda}
\end{align}
and the initial condition
\Beq
  y(0) = \yz.
  \label{cauchy}
\Eeq
\Accorpa\State prima cauchy
\Accorpa\Pbl regy cauchy
Of course, \accorpa{prima}{seconda} can be equivalently rewritten as follows
\begin{align}
  & \ioT \< \dt y , v  > 
  + \intQ \nabla w \cdot \nabla v = 0
  \quad \hbox{for every $v\in\L2 V$},
  \label{intprima}
  \\
  & \intQ wv
  = \intS \dt\yG \, \vG
  + \intQ \nabla y \cdot \nabla v
  + \intS \nablaG\yG \cdot \nablaG\vG
  \qquad
  \non
  \\
  & \quad {}
  + \intQ f'(y) \, v
  + \intS \bigl( \fG'(\yG) - \uG \bigr) \, \vG
  \quad \hbox{for every $(v,\vG) \in\L2 \calV$}.
  \label{intseconda}
\end{align}
\Accorpa\IntPbl intprima intseconda

We are now in a position to introduce our results.
As far as the existence, the uniqueness, the regularity and the continuous dependence 
results are concerned,
we can account for Theorems~2.2, 2.3, 2.4, and 2.6 of \cite{CGS}.
Hence, we have the following statement.
\Bthm
\label{daCGS}
Assume that \HPstruttura, \eqref{hpyz} are fulfilled and let $\uG\in\H1\HG$.
Then, system \Pbl\ admits a unique solution $(y,\yG,w)$ which satisfies
\begin{align}
  & \norma y_{\W{1,\infty}\Vp \cap \H1V \cap \L\infty\Hdue}
  + \norma\yG_{\W{1,\infty}\HG \cap \H1\VG \cap \L\infty\HdueG}
  \non
  \\
  & \quad {}
  + \norma w_{\L\infty V {\cap \L2 {\Hx3}}}
  \,\leq \,C_1,
  \label{stab}
\end{align}
from which, accounting for the Sobolev embedding, it also follows that
\Beq
	\label{y_bounded}
	\norma{y}_{L^\infty(Q)}
	+ \norma{\yG}_{L^\infty(\Sigma)}
	\leq C_1,
\Eeq
for a positive constant $C_1$ that depends only on
$\Omega$, $T$, the shape of the nonlinearities $f$ and~$\fG$,
the initial datum~$\yz$,
and on an upper bound for the norm of $\uG$ in $\H1\HG$. 
Moreover, if $u_{\Gamma\!,i}\in\H1\HG$, $i=1,2$, are two forcing terms
and $(y_i,y_{\Gamma\!,i},w_i)$ are the corresponding solutions,
we have that
\begin{align}
  & \norma{y_1-y_2}_{\L\infty\Vp}^2
  + \norma{y_{\Gamma\!,1}-y_{\Gamma\!,2}}_{\L\infty\HG}^2
    + \norma{\nabla(y_1-y_2)}_{\L2H}^2
	\non
  \\
  & \quad {}
  + \norma{\nablaG(y_{\Gamma\!,1}-y_{\Gamma\!,2)}}_{\L2\HG}^2
  \leq 
  C_2
  \, \norma{u_{\Gamma\!,1}-u_{\Gamma\!,2}}_{\L2\HG}^2,
  \label{contdip}
\end{align}
where the constant $C_2$ depends only on
$\Omega$, $T$, and the shape of the nonlinearities $f$ and~$\fG$.
\Ethm

Once the well-posedness of the system \Pbl~has been proved, we can
address the corresponding control problem. 
As \andrea{far as} the assumptions on the cost functional are concerned,
we postulate that
\begin{align}
  &
  \label{hpzzzz}
  \zQ \in \H1 H , \ \ 
  \zS \in \LS2 , \ \ 
  \zO \in \Hx1, \ \ 
  \zG \in \HxG 1.
  \\ 
   \label{hp_bb}
	& \bQ,\, \bS,\, \bO,\, \bG,\, \bz
	\ \ \hbox{are nonnegative constants, but not all zero.} 
	\\  \non
	& \Mz > 0 , \ \ 
	\uGmin ,\, \uGmax \in \LS\infty ,
	\ \ \hbox{with } \ \
	\uGmin \leq \uGmax \quad \andrea{\aeS},
	\\
	& \quad \hbox{in such a way that $\Uad$ turns out to be nonempty.}
	\label{hpUad}
\end{align}
Below, the first fundamental result related to the existence of optimal controls
\andrea{can be found}.
\Bthm
\label{Optimum}
Assume that \HPstruttura, \eqref{hpyz}, and \accorpa{hpzzzz}{hpUad} are in force.
Then, there exists $\uopt\in\Uad$ such~that
\Beq
  \calJ(\yopt,\yGopt,\uopt)
  \leq \calJ(y,\yG,\uG)
  \quad \hbox{for every $\uG\in\Uad$},
  \label{optimum}
\Eeq
where $\yopt$, $\yGopt$ and $y$, $\yG$
are the components of the solutions $(\yopt,\yGopt,\wopt)$ and $(y,\yG,w)$
to~the state system \Pbl\ corresponding to the controls
$\uopt$ and~$\uG$, respectively.
Such a control variable $\uopt$ is called optimal control.
\Ethm
%
%
The well-posedness of the system \Pbl, 
allows us to properly define the \socal~control-to-state mapping. We set
\begin{align}
  \non
  & \calX := \H1\HG \cap \LS\infty
  \ \ \hbox{and} \  \
  \calY := \H1\calG \cap \L\infty\calV,
  \non
  \\
  & \hbox{$\calU$ is an open bounded set in $\calX$ that includes $\Uad$,}
  \non
  \\
  & \calS : \calU \subset \calX \to \calY \, \hbox{ is defined by } \calS(\uG) := (y,\yG), 
  \non
  \\  \non
  &  \quad \hbox{where }\,\hbox{$(y,\yG,w)$ is the}
  \hbox{ solution to \Pbl\ corresponding to }\uG.
\end{align}
\Brem
Note that the existence of the superset $\cal U$ containing $\Uad$ is trivially satisfied. 
Indeed, for instance, we can take
\Beq
	\non
	{\cal U}:=\bigl\{ \uG \in {\cal X} : \norma{\uG}_{L^\infty (\Sigma)} 
	< \norma{\uGmin}_{L^\infty (\Sigma)}+\norma{\uGmax}_{L^\infty (\Sigma)}+1, \ 
	\norma{\dt\uG}_{L^2(\Sigma)}< M_0+1 \bigr\}.
\Eeq
\Erem

Thus, we can express the cost functional $\cal J$ as a function of $\uG$ by introducing
the \socal~reduced cost functional
\Beq
  \redJ : \calU \to \erre
  \quad \hbox{ which is defined by} \quad
  \redJ(\uG) := \calJ(\calS(\uG),\uG).
  \label{defredJ}
\Eeq
Formally, as $\Uad$ is convex, it is a standard matter to realize that 
the desired necessary condition for $\uopt$ is carried out by the following
variational inequality
\Beq
  \< D\redJ(\uopt) , \vG-\uopt > \geq 0
  \quad \hbox{for every $\vG\in\Uad$},
  \label{precnopt}
\Eeq
where $D\redJ(\uopt)$ denotes the derivative of $\redJ $ at 
$\uopt$ in a suitable functional sense.
The strategy we follow in order to obtain some optimality conditions consists 
in proving at first that $\calS$ is \Frechet\ differentiable at $\uopt$, \andrea{and}
then, accounting for the chain rule, developing the above inequality
to get an explicit formulation which characterizes the optimality.
As we shall see in Section \ref{FRECHET}, this procedure naturally leads to 
the linearized system, that we briefly introduce in the lines below. 
Let us fix $\uopt\in\calU$, the corresponding state $(\yopt,\yGopt):=\calS(\uopt)$,
and introduce the increment $\hG\in\H1\HG$.
\andrea{Moreover}, we set for convenience
\Beq
  \lam := f''(\yopt),
  \aand
  \lamG := \fG''(\yGopt).
  \label{deflam}
\Eeq
\andrea{Then}, the linearized system for \Ipbl~consists of finding a triplet $(\xi,\xiG,\eta)$ 
satisfying the analogue of \Regsoluz, solving\andrea{,} \aat\andrea{,} the variational equations
\begin{align}
  &\< \dt\xi(t) , v >
  + \iO \nabla\eta(t) \cdot \nabla v = 0
  \quad \hbox{for every $v\in V$},
  \label{linprima}
  \\
    & \iO \eta(t) v
  = \iG \dt\xiG(t) \, \vG
  + \iO \nabla\xi(t) \cdot \nabla v
  + \iG \nablaG\xiG(t) \cdot \nablaG\vG
  \non
  \\
  & \quad
  + \iO \lam(t) \, \xi(t) \, v
  + \iG \bigl( \lamG(t) \, \xiG(t) - \hG(t) \bigr) \, \vG
  \quad \hbox{for every $(v,\vG) \in\calV$}
  \label{linseconda}
\end{align}
and satisfying the initial condition
\Beq
  \xi(0) = 0 .
  \label{lincauchy}
\Eeq
\Accorpa\Linpbl linprima lincauchy
In order to obtain the well-posedness for the above system, we 
would be tempted to directly invoke \cite[Thm.~2.4, p.~978]{CGS}. However, from
a careful investigation, we realized that the requirements on $\lam$
are not satisfied in our setting. Indeed, note that $\dt \lam = f'''(\yopt)\dt\yopt$
and that $\dt \yopt\in \L\infty \Vp \cap \L2 V$ \andrea{by virtue of Theorem~\ref{daCGS}}. 
So, in our framework we cannot infer that $\lam \in W^{1,\infty}(0,T;H)$.
This lack of regularity, due to the absence of the viscous term, can
be however overcome by applying a different estimate in the term involving 
$\lam$. 
Therefore, modifying properly the proof of \cite[Thm.~2.4, p.~978]{CGS},
the same result holds. 
\Bthm
\label{Existlin}
Let $\uopt\in\calU$, $(\yopt,\yGopt)=\calS(\uopt)$, 
and $\lam$, $\lamG$ be defined by~\eqref{deflam}.
Then, for every $\hG\in\H1\HG$,
there exists a unique triplet $(\xi,\xiG,\eta)$
satisfying the analogue of \Regsoluz\
and solving the linearized system \Linpbl.
\Ethm

Next, we will show that $\calS$ is \Frechet\ differentiable at $\uopt$,
that $D\calS(\uopt)$ is a linear operator from $\calX$ to $\calY$,
and also that, for every $\hG$ in $\calX$,
$
	[D\calS(\uopt)](\hG)=(\xi,\xiG),
	$
where the triplet $(\xi,\xiG,\eta)$ represents the unique solution 
to the linearized system associated to $\hG$.
Here is the precise result.
\Bthm
\label{Fdiff}
Let $\uopt\in\calU$, $(\yopt,\yGopt)=\calS(\uopt)$, and $\lam$\andrea{,} $\lamG$ be defined by~\eqref{deflam}.
Then the control-to-state mapping $\calS:{\cal U}\subset {\cal X}\to{\cal Y} $ is \Frechet\ 
differentiable at $\uopt$. Moreover, its derivative 
$D\calS(\uopt)$ is a linear operator from $\cal U$ to $\calY$ which is given as follows:
whenever $\hG\in\calX$ fulfills $\uopt +\hG \in \cal U$, 
the value of $D\calS(\uopt)$ at $\hG$ consists of the pair $(\xi,\xiG)$,
where $(\xi,\xiG,\eta)$ is the unique solution to the linearized system \Linpbl.
\Ethm
Then, by invoking the chain rule, we develop \eqref{precnopt} in order to 
obtain the following explicit optimality condition
\begin{align}
  & \bQ \intQ (\yopt - \zQ) \xi
  + \bS \intS (\yGopt - \zS) \xiG
  + \bO \iO (\yopt(T) - \zO) \xi(T)
  \non
  \\
  & \quad {}
  + \bG \iG (\yGopt(T) - \zG) \xiG(T)
  + \bz \intS \uopt (\vG - \uopt)
  \geq 0
  \quad \hbox{for every $\vG\in\Uad$},
  \qquad
  \label{cnopt}
\end{align}
where $\xi$ and $\xiG$ are the first two components
of the unique solution to the linearized system 
corresponding to $\hG=\vG-\uopt$.

Lastly, we try to eliminate the pair $(\xi,\xiG)$ from the above inequality.
To overcome this issue, we introduce the \socal~adjoint system.
\andrea{Namely}, we are looking for a triplet $(q,\qG,p)$ that fulfills the regularity requirements
\begin{align}
   (q,\qG ) &\in \H1 \GO \cap \L\infty \VO \cap \L2 {\Hx2 \times H^2(\Gamma)},
  \label{qqG}
  \\
  p &\in \H1 V \dafare{\cap \L2 {\Hx4},}
  \label{regp}
  \\    
  \qG(t) &= q(t)\suG  \,\, \aat\andrea{,}
\end{align}
\Accorpa\Regadj qqG regp
\pier{and solves}, $\aat$, the following backward-in-time problem
\begin{align}
  \label{primaadjpure}
  & \iO q(t) \, v
  = \iO \nabla p(t) \cdot \nabla v
  \quad \hbox{for every $v\in V,$}
  \\
  & - \iO \dt p(t) \, v
  + \iO \nabla q(t) \cdot \nabla v
  + \iO  \lam(t) q(t)  v
  - \iG \dt \qG(t) \, \vG
  + \iG \nablaG  \qG(t) \cdot \nablaG\vG
  \non
  \\ \non
  &  \quad\quad\quad\quad
  + \iG  \lamG(t) \qG(t) \vG
  = \bQ \iO (\yopt(t) - \zQ(t)) v
  + \bS \iG (\yGopt(t) - \zS(t)) \vG
  \\
  & \hspace{9.8cm} \hbox{for every $(v,\vG) \in\calV$},
  \label{secondaadjpure}
\end{align}
and the final condition
\begin{align}
   \non & \iO p(T) v 
   + \iG \qG(T) \vG
   = \bO \iO (\yopt(T) - \zO) v(T) 
   + \bG \iG (\yGopt(T) - \zG) \vG(T)
   \\
   &
    \hspace{9.6cm} \hbox{for every $(v,\vG) \in\calV$}.
   \label{terzaadjpure}
\end{align}
\Accorpa\Pbladj primaadjpure terzaadjpure
In order to simplify the notation, let us convey to denote
\begin{align*}
	\phi_Q:=\bQ (\yopt - \zQ), \ \ 
	\phi_\Sigma:=\bS (\yGopt - \zS), \ \
	\phi_\Omega:=\bO (\yopt(T) - \zO), \ \
	\phi_\Gamma:=\bG (\yGopt(T) - \zG).
\end{align*}
Here, the well-posedness result follows.
\Bthm
\label{adjexist}
Let $\uopt$ be an optimal control with the corresponding optimal state
$(\yopt,\yGopt)$.
Moreover, let us postulate that $\phO$ and $\phG$ satisfy the following compatibility condition:
there exists a couple $(\Phi, \phi_\Gamma)\in \VO$
such that $\phi_\Omega={\cal N }(\Phi) + (\phi_\Omega)^\Omega$.
Then, the adjoint system \Pbladj~admits a unique solution $(p,q,\qG)$
satisfying the regularity requirements \Regadj.
\Ethm

\andrea{Let us underline that the above result is new with respect
to \cite{CGSpure}, where just the existence, in a very weak setting, 
was proved. 
Here, the complete well-posedness of the adjoint system
is now achievable under the enforced assumptions on the potential setting.
Moreover, notice that the unique solution to \Pbladj\ enjoys the strong 
regularity \Regadj. 
Then, once}
that the adjoint variables are at our disposal, we are in a position
to eliminate $\xi$ and $\xiG$ from \eqref{cnopt}, thus leading to the following
optimality condition.
\Bthm
\label{CNoptadj}
Let $\uopt$ be an optimal control, $(\yopt,\yGopt)$ be the corresponding 
optimal state, and $(p,q,\qG)$
be the associated solution to the adjoint system \accorpa{regp}{secondaadjpure}. 
Then, the first-order necessary condition for optimality is characterized by the
following variational inequality
\Beq
  \intS (\qG + \bz \uopt) (\vG - \uopt) \geq 0
  \quad \hbox{for every $\vG\in\Uad$}.
  \label{cnoptadj}
\Eeq
Moreover, whenever $\bz>0$, it turns out that
\Beq
  \hbox{\sl\mathsurround 3pt
  $\uopt$ is the orthogonal projection of $-\qG/\bz$ on $\Uad$}
  \label{projection}
\Eeq
with respect to the standard inner product of~$\LS2$.
\Ethm

\pier{%
\Brem
Of course, the condition \eqref{cnoptadj} 
also \lastrev{entails} that the element 
$ - (\qG + \bz \uopt)$ belongs to the {\it normal cone} 
of the closed and convex set $\Uad$ (defined in~\eqref{defUad}) at $\uopt$
in the framework of the Hilbert space~$\LS2$.
Owing to the structure of the control-box $\Uad$, 
if $\bz>0$ then the projection $\uopt$ in \eqref{projection}
is the one among the elements $\uG\in\H1 {L^2(\Gamma)}\cap\LS\infty$ satisfying the two constraints 
%
$$ \uGmin\leq\uG\leq\uGmax\quad \aeS,\ \quad \norma{\dt\uG}_{L^2(\Sigma)}\leq\Mz $$
that is closest 
to $-\qG/\bz$ in the sense of the norm in~$\LS2$. In particular, 
if the function $z_\Gamma $ defined by 
\begin{align}
	\non
	z_\Gamma (x,t)=\max \bigl\{ 
	\uGmin(x,t), \min\graffe{\uGmax(x,t), - \qG(x,t)/\bz}
	\bigr\} , \quad (x,t) \in \Sigma ,
\end{align}
belongs to $\H1{L^2(\Gamma)}$ and its time derivative fulfills $\norma{\dt z_\Gamma}_{L^2(\Sigma)}\leq\Mz$, then we necessarily have that $\lastrev{\uopt} = z_\Gamma $ \ $\aeS$. 
\Erem}

In the remainder, we introduce further notation and recall some 
\wk~inequalities and general facts which will be useful later on.
First of all, we often owe to the Young inequality
\Beq
  ab \leq \delta a^2 + \frac 1 {4\delta} \, b^2
  \quad \hbox{for every $a,b\geq 0$ and $\delta>0$}.
  \label{young}
\Eeq
Furthermore, we account for the Poincar\'e inequality
\Beq
  \normaV v^2 \leq C_{\Omega} \bigl( \normaH{\nabla v}^2 + |\vO|^2 \bigr)
  \quad \hbox{for every $v\in V$},
  \label{poincare}
\Eeq
where $C_{\Omega}$~depends only on~$\Omega$.
Furthermore, we point out the following inequality, to which we will 
refer to as compactness inequality (see, e.g., \cite[Lem.~5.1, p.~58]{LionsQM}):
for every $\delta >0$ there exists $\cd>0$ such that
\Beq
	\label{compactineq}
	\normaH{v}^2
	\leq
	\delta \norma{\nabla v}_{H}^2
	+
	c_\delta \norma{v}_{\Vp}^2
	\quad \hbox{for every $v\in V$},
\Eeq
where the constant $c_\delta$ depends only on $\delta$ and $\Omega$.

Lastly, let us point out a convention we use in the 
whole paper as far as the constants are concerned.
We agree that the small-case symbol $c$ stands for different constants depending only
on the final time~$T$, on~$\Omega$, the shape of the nonlinearities
and the norms of functions involved in the assumptions of our statements. 
For this reason, its meaning might
change from line to line and even in the same chain of calculations.
Conversely, the capital letters are devoted to denote precise constants 
which we eventually will refer to.

\section{Existence of an optimal control}
\label{OPTIMUM}
\setcounter{equation}{0}
From this section on, we \andrea{will} start with the proofs of the stated results.
Here, we aim to prove the existence of optimal control.
Before moving on, let us briefly remark that Theorem \ref{daCGS} 
is slightly stronger with respect to the result of \cite{CGS}, 
since by \eqref{regw}
we require additional space regularity for the variable $w$.
As a matter of fact, it suffices to combine the original result with a comparison
argument to realize that $w\in \L2 {\Hx3}$ as well.
\Bdim[Proof of Theorem~\ref{daCGS}]
Since the proof is the same as in \cite{CGS},
we can afford to be sketchy by just pointing out some highlights.
From \cite{CGS}, it follows that there exists a positive constant $c$ such that
\begin{align*}
  & \norma y_{\W{1,\infty}\Vp \cap \H1V \cap \L\infty\Hdue}
  + \norma\yG_{\W{1,\infty}\HG \cap \H1\VG \cap \L\infty\HdueG}
  \non
  \\
  & \quad {}
  + \norma w_{\L\infty V }
  \,\leq c.
\end{align*}
On the other hand, owing to the above estimate,
by comparison in equation \eqref{Iprima}, we infer that
$\Delta w  \in \L2 V$, and the classical elliptic regularity theory
directly implies that $w\in \L2 {\Hx3}$.
\Edim
\proof[Proof of Theorem~\ref{Optimum}]
We proceed by employing the direct method. First, let us pick a minimizing sequence
$\graffe{\unG}_n$ for the cost functional $\cal J$ and, for every $n$, 
let us denote by $(\yn,\ynG,\wn)$ the corresponding solution to \Pbl.
Since, for every $n$, $\unG$ belongs to $\Uad$ and 
the triplet $(\yn,\ynG,\wn)$ solves the state system, 
the bounds \accorpa{stab}{y_bounded} are in force.
Thus, for every $n$, we infer that 
\Beq
	\rmin  \leq \yn \leq  \rmax
	\quad \hbox{\aeQ},
	\quad
	\rmin  \leq y_{\Gamma,n} \leq  \rmax
	\quad \hbox{\aeS},
  \label{dafaraway}
\Eeq
for some $\rmin, \rmax$ \andrea{satisfying} $-\infty < \rmin \leq \rmax < + \infty$.
It is now a standard matter to show that, accounting for weak and weak-star
compactness arguments (see, e.g., \cite[Sect.~8, Cor.~4, p. 85]{Simon}),
up to a subsequence, the following convergences are verified
\Bsist
  & \unG \to \uopt
  & \quad \hbox{weakly star in $\LS\infty\cap\H1H$},
  \non
  \\
  & \yn \to \yopt
  & \quad \hbox{weakly star in $\W{1,\infty}\Vp \cap \H1V \cap \L\infty\Hdue$}
  \non
  \\
  && \qquad \hbox{and strongly in $\C0V$},
  \non
  \\
  & \ynG \to \yGopt
  & \quad \hbox{weakly star in $\W{1,\infty}\HG \cap \H1\VG \cap \L\infty\HdueG$}
  \non
  \\
  && \qquad \hbox{and strongly in $\C0\VG$},
  \non
  \\
  & \wn \to \wopt
  & \quad \hbox{weakly star in $\L\infty V\cap \L2 {\Hx3}$}.
  \non
\Esist
In addition, since $\Uad$ is closed \andrea{it follows} that $\uopt\in\Uad$ and,
from the strong convergences pointed out above that $\yopt(0)=\yz$.
Moreover, the strong convergences of $\yn$ and $\ynG$, combined with
the regularity of $f$ and $\fG$, imply that
\begin{align}
	 \non
	f'(\yn) &\to  f'(\yopt) \quad \hbox{strongly in $\C0H$},
	\\  \non
	\fG'(\ynG) &\to  \fG'(\yGopt) \quad \hbox{strongly in $\C0 \HG$}.
\end{align}
By virtue of all these convergences, we can easily pass to the limit in the 
integrated variational formulation~\IntPbl\
written for $(\yn,\ynG,\wn)$ and $\unG$ 
to conclude that $(\yopt,\yGopt,\wopt)$ 
solves \IntPbl\ with $\uG:=\uopt$.
Lastly, \andrea{accounting for} the lower weak semicontinuity of $\cal J$,
is \sfw~to realize that $(\yopt,\yGopt,\uopt)$ 
is \andrea{indeed} the minimizer we are looking~for.
\qed

\section{The control-to-state mapping}
\label{FRECHET}
\setcounter{equation}{0}
In this section, we first prove Theorem \ref{Existlin}, and then
show the \Frechet~differentiability of the control-to-state 
\andrea{operator} $\calS$ between suitable Banach spaces.

\begin{proof}[Proof of Theorem~\ref{Existlin}]
As sketched above, we would like to \andrea{invoke}
\cite[Thm.~2.4, p.~978]{CGS}. Unfortunately, the assumption $\lam\in W^{1,\infty}(0,T;H)$ fails
to be satisfied. On the other hand, due to the regularity of the potential $f$, 
along with \eqref{y_bounded}, we can easily check that
$\dt \lam = f'''(\yopt)\dt\yopt$ belongs at least to $\L2 H$. 
Let us claim that this regularity is
actually sufficient in order to prove the same result as in \cite[Thm.~2.4]{CGS}.
Since it consists of a minor change, let us proceed quite formally, 
leaving the details to the reader and avoiding to write the explicit dependence 
on the time variable for convenience.
As a starting point let us assume that
\begin{align*}
	(\xi,\xiG) &\in \H1 {{\cal G}} \cap \L\infty {{\cal V}} \cap \L2 {\Hx2 \times \HxG2},
	\\
	\eta & \in \L2 V,
\end{align*}
which can be easily obtained applying \cite[Thms.~2.2 and 2.3]{CGS}.
Moreover, note that the mean value of $\xi$ is \andrea{zero}, thanks to 
\eqref{linprima} and \eqref{lincauchy}.
Then, let us formally differentiate 
\accorpa{linprima}{linseconda} with respect to time and integrate over $(0,t)$ to get
\begin{align*}
  & 
  	\iot \< \dtt \xi, v>
	+ \intQt \nabla \dt \eta \cdot \nabla v
	= 0
  \quad \hbox{for every $v\in V$},
  \\
  & \intQt \dt\eta v
  = \intSt \dtt\xiG \, \vG
  + \intQt \nabla\dt\xi \cdot \nabla v
  + \intSt \nablaG\dt\xiG \cdot \nablaG\vG
  + \intQt \lam \, \dt\xi \, v
  \non
  \\
  & \quad
  + \intQt \dt\lam \, \xi \, v
  + \intSt \lamG \dt\xiG \vG
  + \intSt \dt\lamG  \xiG \vG
  - \intSt \dt\hG  \vG
 \quad 
 \hbox{for every $(v,\vG) \in\calV$}.
\end{align*}
Next, we test the former by ${\cal N}(\dt \xi)$ and the latter by $-(\dt\xi,\dt\xiG)$,
add the resulting equalities and integrate by parts to obtain, after some simplifications,
that
\begin{align}
  	\non
  	& 
  	\frac 12 \norma{\dt\xi(t)}_*^2
	+ \frac 12 \iG |\dt\xi(t)|^2
	+ \intQt |\nabla \dt \xi|^2
	+ \intSt |\nablaG \dt \xiG|^2
	\\
	&\quad \non
	=  
	\frac 12 \norma{\dt\xi(0)}_*^2
  	+ \frac 12 \iG |\dt\xi(0)|^2
	- \intQt \lam |\dt\xi|^2
  - \intSt \lamG |\dt\xiG|^2
  \\
	&\qquad
  - \intQt \dt\lam \, \xi \, \dt\xi
  - \intSt \dt\lamG \, \xiG \, \dt\xiG
  + \intSt \dt\hG  \dt\xiG.
  \label{pier1}
\end{align}
Let us denote the terms on the \rhs~by $I_1,...,I_7$, in this order.
Owing to the Young inequality and to 
the boundedness of $\lamG$, we easily handle the boundary terms
as follows
\begin{align*}
	|I_4|+|I_7|
	\leq
	c \intSt |\xiG|^2
	+ c \intSt |\dt\xiG|^2
	+ c \intSt |\dt\hG|^2,
\end{align*}
where let us remark that $\xiG$ has been already estimated in $\H1 {\HG}$.
Moreover, owing to the inequality \eqref{compactineq} and to the \Poincare~inequality
\eqref{poincare}, we have that
\begin{align*}
	|I_3|
	\leq
	c \intQt |\dt\xi|^2
	\leq 
	\frac 14 \intQt |\nabla \dt \xi|^2
	+ c \iot \norma{\dt\xi}_*^2.
\end{align*}
Finally, using the \Holder\ inequality and \eqref{compactineq}, we obtain that
\begin{align*}
	|I_5|
	&\leq
	c \iot \norma{\dt\lam}_2 \norma{\xi}_4 \norma{\dt\xi}_4
	\leq 
	c 	\iot \norma{\dt \lam}_H\norma{\xi}_V \norma{\dt\xi}_V
	\\ 
	& 
	\leq
	\frac 18 \iot \norma{\dt\xi}_V^2
	+ c 	\iot \norma{\dt \lam}_H^2\norma{\xi}_V^2
	\leq 
	\frac 18 \intQt |\dt\xi|^2
	+ \frac 18 \intQt |\nabla\dt \xi|^2
	+ c \intQt |\dt \lam|^2
	\\ &
	\leq 
	\frac 14 \intQt |\nabla\dt \xi|^2
	+ c \iot \norma{\dt\xi}_*^2
	+ c,
\end{align*}
thanks to the Sobolev embedding
$V \subset \Lx4$ and the fact that $\xi$ has already been estimated in 
$\L\infty V$.
\dafare{
In a similar way, we can deal with $I_6$ by using the \Holder\ and
Young inequalities. In fact, we have that 
\begin{align*}
	&
	|I_6|
	\leq
	c \iot \norma{\dt\lamG}_{\LxG 2} \norma{\xiG}_{\LxG 4}\norma{\dt\xiG}_{\LxG 4}
	\leq
	c \iot \norma{\dt\lamG}_{\HG} \norma{\xiG}_{\VG}\norma{\dt\xiG}_{\VG}
	\\ & \quad
	\leq
	\frac 12 \iot \norma{\dt\xiG}_{\VG}^2
	+ c \iot \norma{\xiG}_{\VG}^2\norma{\dt\lamG}_{\HG}^2
	\\ & \quad
	\leq 
	\frac 12 \intSt | \nabla \dt\xiG|^2
	+ \frac 12 \intSt | \dt\xiG|^2
	+ c \norma{\xiG}_{\L\infty \VG}^2\norma{\dt \lamG}_{\L\infty \HG}^2
	\\ & \quad
	\leq 
	\frac 12 \intSt | \nabla \dt\xiG|^2
	+ c,
\end{align*}
where we apply the Sobolev embedding $\VG \subset \LxG 4$, the fact that
$\dt\lamG \in \L\infty \HG$, and that $\xiG$ has
been already estimated in $\H1 \HG \cap \L\infty \VG$.}
Therefore, it suffices to show that $I_1$ and $I_2$ remained bounded.
In this regards, we evaluate equations \eqref{linprima} and \eqref{linseconda} 
at $t=0$. Then, we test them by ${\cal N}(\dt\xi(0))$ and $-\dt\xi(0)$, add the 
resulting equalities and rearrange the terms to obtain that
\begin{align*}
	\norma{\dt\xi(0)}_*^2
	+ \norma{\dt\xiG(0)}_{L^2(\Gamma)}^2
	= \iG \hG (0) \dt\xiG(0),
\end{align*}
where the initial condition $\xi(0)=0$ has been exploited.
Hence, we use the Young inequality to handle the term on the 
\rhs~and infer that
\begin{align*}
	\norma{\dt\xi(0)}_*^2
	+ \frac 12 \norma{\dt\xiG(0)}_{L^2(\Gamma)}^2
	\leq 
	\frac 12 \norma{\hG(0)}_{\HG}^2
	\leq 
	c \norma{\hG}_{\C0 {\HG}}^2
	\leq
	c \norma{\hG}_{\H1 {\HG}}^2
	\leq 
	c.
\end{align*}
Then, recalling \eqref{pier1} and collecting
the previous estimates, we conclude the proof by
applying the Gronwall lemma.
\end{proof}

We will see that, in order to directly check the definition of \Frechet~differentiability
for $\cal S$, some \andrea{stronger} continuous dependence results 
\andrea{with respect to \eqref{contdip}}
need to be shown. 
Therefore, this is the task of the following lemmas.
\andrea{The first one is somehow the corresponding non-viscous 
version of \cite[Lem.~4.1, p.~207]{CGS_OPT}.}
\Blem
\label{Stability}
Let $u_{\Gamma,i}\in \calU \,$ for $i=1,2$
and let $(y_i,y_{\Gamma,i},w_i)$ be the corresponding solutions
to \accorpa{regy}{cauchy}.
Then, it follows that
\Beq
  \norma{(y_1,y_{\Gamma,1}) - (y_2,y_{\Gamma,2})}_{\calY}
  \leq 
  C_3 \norma{u_{\Gamma,1} - u_{\Gamma,2}}_{\L2\HG},
  \label{stability}
\Eeq
for a positive constant $C_3$ that may depend on
$\Omega$, $T$, the shape of the nonlinearities $f$ and~$\fG$,
and on the initial datum~$\yz$.
\Elem

\Bdim
To begin with, let us fix for convenience the notation
\Beq
  \uG := u_{\Gamma,1} - u_{\Gamma,2}, \quad
  y := y_1 - y_2 , \quad
  \yG := y_{\Gamma,1} - y_{\Gamma,2},
  \aand
  w := w_1 - w_2.
  \label{differenze}
\Eeq
Then, we write the system \State~for both the solutions $(y_i,y_{\Gamma,i},w_i)$
for $i=1,2$, and take the difference to obtain, \aat, that
\begin{align}
	& \< \dt y(t) , v >
  + \iO \nabla w(t) \cdot \nabla v = 0
  \quad \hbox{for every $v\in V$},
  \label{primadiffdue}
  \\
  & \iO w(t) \, v
  = \iG \dt\yG(t) \, \vG
  + \iO \nabla y(t) \cdot \nabla v
  \non
  \\
  \non
  & \quad {}
  + \iG \nablaG\yG(t) \cdot \nablaG\vG
  + \iO \bigl( f'(y_1(t)) - f'(y_2(t)) \bigr) \, v
  \\
  & \quad{} + \iG \bigl( \fG'(y_{\Gamma,1}(t)) - \fG'(y_{\Gamma,2}(t)) - \uG(t) \bigr) \, \vG
    \quad \hbox{for every $(v,\vG) \in\calV$}
	\label{secondadiff}
\end{align}
and $y(0)=0$.
Moreover, we point out that $\dt y$ has zero mean value 
since \eqref{conserved} holds for both $\dt y_1$ and $\dt y_2$
so that $\calN(\dt y)$ can be considered as a test function.
So, we subtract to both \andrea{sides} of \eqref{secondadiff} the terms $\iO y(t) \, v $ and 
$\iG \yG(t) \, \vG $, write the above equations at the time $s$, 
test \eqref{primadiffdue} by $\calN(\dt y(s))$, the new \eqref{secondadiff} 
by $-\dt(y,\yG)(s)$, add the resulting equalities and integrate over $(0,t)$
for an arbitrary $t\in(0,T)$. We obtain that
\begin{align}
  \non
  & \iot \< \dt y , \calN(\dt y) > 
  + \intQt \nabla w \cdot \nabla\calN(\dt y)
  - \iot \< \dt y , w > 
  + \intSt |\dt\yG|^2
  \\
  \non
  & \qquad {} 
  + \frac 12 \normaV{y(t)}^2
  + \frac 12 \normaVG{\yG(t)}^2
   = - \iot \< \dt y,  f'(y_1) - f'(y_2) - y > 
     \\
  & \qquad {} 
  - \intSt \bigl( \fG'(y_{\Gamma,1}) 
  - \fG'(y_{\Gamma,2}) - \yG \bigr) \dt\yG 
  + \intSt \uG \, \dt\yG,
  \label{perstab}
\end{align}
where we also invoke the fact that $y(0)=0$, $\yG(0)=0$ since $y_1$ and $y_2$
have the same initial value $y_0$.
The first three integrals of the above equality can be treated with
the help of \eqref{defN} and \eqref{simmN} as follows
\begin{align}
  \non
  \iot \< \dt y , \calN(\dt y) >
  + \intQt \nabla w \cdot \nabla\calN(\dt y)
  -  \iot \< \dt y , w > 
  = \iot \norma{\dt y}_*^2 \geq \,0 \,.
\end{align}
Furthermore, all the other contributions on the \lhs~are nonnegative,
so we \andrea{are reduced} to control the integrals on the \rhs. 
On the other hand, both $y_1$ and $y_2$, as solutions to \Ipbl, satisfy~\eqref{stab} 
and \eqref{y_bounded}.
Using the Young inequality, we estimate the first term 
of the \rhs~by
\begin{align}
  - \iot \<  \dt y , f'(y_1) - f'(y_2) - y > 
  \non
   &\leq \iot \normaVp{\dt y} \normaV{f'(y_1) - f'(y_2) - y}
  \\
  \label{primoint}
  &  
  \leq \frac 12 \iot \normaVp{\dt y}^2 
  + \frac 12 \iot \normaV{f'(y_1) - f'(y_2) - y}^2.
\end{align}
By invoking the \Lip~continuity of $f'$ and $f''$, and the Sobolev embedding $V\subset \Lx4$,
we are able to bound the last term of the previous estimate as follows
\begin{align}
  \normaV{f'(y_1) - f'(y_2) -y }^2 
  & \leq c\normaV y ^2 + c\normaV{f'(y_1) - f'(y_2) }^2 
  \non
  \\
  &  \leq c\normaV y ^2 + c\normaH{f'(y_1) - f'(y_2)}^2 
  +c\normaH{\nabla \bigl(f'(y_1) - f'(y_2) \bigr)}^2
  \non
  \\
  &  \leq c\normaV y ^2 + c\normaH y ^2 
  + \bigl( \normaH {f''(y_1)\nabla y}^2 
  + \normaH{\bigl( f''(y_1) - f''(y_2)\bigr) \nabla y_2}^2 \bigr)
  \non
  \\ \non
  &  \leq c\normaV y ^2 + c \normaH{\nabla y} ^2 
  + c \iO |y|^2 |\nabla y_2|^2
  \\ \non
  & 
  \leq 
  c \norma{y}_V^2 + c \norma{y}^2_4\norma{\nabla y_2}^2_4
  \leq c \normaV{y}^2 + c \normaV{y}^2 \norma{y_2}_{\Hdue}^2  ,
\end{align}
where in the last inequality we invoke the fact that, as a solution,
$y_2$ satisfies \eqref{stab} so that $y_2$ is bounded in $\L\infty \Hdue$.
Summing up, the estimate
\Bsist
  \non
  && - \iot \< \dt y , f'(y_1) - f'(y_2) - y > 
  \leq \frac 12 \iot \normaVp{\dt y}^2
  + c \iot  \normaV{y}^2
\Esist
has been shown.
\andrea{The boundary integrals can be easily handled
owing to \eqref{young} and the
\Lip~continuity of $\fG$. Indeed, we have that}
\Beq
  - \intSt \bigl( \fG'(y_{\Gamma,1}) - \fG'(y_{\Gamma,2}) - \yG \bigr) \dt \yG
  \leq \frac 14 \intSt |\dt \yG|^2
  + c \intSt |\yG|^2,
  \non
\Eeq
and
\Beq
  \intSt \uG \, \dt\yG
  \leq \frac 14 \intSt |\dt\yG|^2
  + c \intSt |\uG|^2,
  \non
\Eeq
respectively. Lastly, upon collecting all the previous estimates, we realize that
\begin{align}
  & \frac 12 \iot \normaVp{\dt y}^2
  + \frac 12 \intSt |\dt\yG|^2
  + \frac 12 \normaV{y(t)}^2
  + \frac 12 \normaVG{\yG(t)}^2
  \non
  \\
  & \quad{} \non
  \leq c \iot \normaV{y}^2
  + c \intSt  |\yG|^2
  + c \intSt |\uG|^2 ,
\end{align}
whence the standard Gronwall lemma yields the stability inequality we are looking for.
\Edim
Unfortunately, we will see that in order to prove the \Frechet~differentiability of $\calS$,
the above \pier{lemma} turns out to be insufficient. Then, 
in the lines below we present an improvement.
\andrea{Notice that the following result is a novelty in comparison
to \cite{CGS_OPT}. On the other hand, it turns out to be necessary
in order to handle the control problem we are dealing with.
}
\Blem
\label{Stability2}
Let $u_{\Gamma,i}\in \calU$ for $i=1,2$
and let $(y_i,y_{\Gamma,i},w_i)$ be the corresponding solutions
to \Ipbl. 
Then, there exists a positive constant $C_4$ such that
\begin{align}
  & \non
  \norma{(y_1,y_{\Gamma,1}) - (y_2,y_{\Gamma,2})}_{\W{1,\infty}\calG 
  \cap \H1\calV \cap \L\infty{\Hdue \times \HdueG}}
  \\ & \quad 
	\leq 
  C_4
  \norma{u_{\Gamma,1} - u_{\Gamma,2}}_{\H1\HG},
  \label{stability2}
\end{align}
where $C_4$ is a positive constant which depends only on
$\Omega$, $T$, the shape of the nonlinearities $f$ and~$\fG$, 
the initial datum~$\yz$.
\Elem
\Bdim
In what follows, to keep the proof as easy as possible, we proceed formally.
The justification can be carried out rigorously, e.g.,
within a time-discretization scheme.
Then, providing to show some estimates for the differences, one has to pass 
to the limit in suitable topologies. 

To begin with, we write the problem \State~for both the solutions 
$(y_i,y_{\Gamma,i},w_i)$, $i=1,2$, take the difference 
and use the notation set by \eqref{differenze}.
Then, we differentiate the equations with respect to the time variable to obtain
that, \aat, the following are satisfied
\begin{align}
	\label{abst_var_uno}
  & \< \dtt y(t) , v >
  + \iO \nabla \dt w(t) \cdot \nabla v = 0
  \quad \hbox{for every $v\in V,$ }
  \\
  & \iO \dt w(t) \, v
  = \iG \dtt \yG(t) \, \vG
  + \iO \nabla \dt y(t) \cdot \nabla v
  + \iG \nablaG \dt \yG(t) \cdot \nablaG\vG
  \quad
  \non
  \\
  & \quad {}
  \non
  + \iO f''(y_1(t)) \dt y(t) \, v
  + \iO \bigl( f''(y_1(t)) - f''(y_2(t)) \bigr)\dt y_2(t) \, v
  \\
 & \quad
  \non
  + \iG \fG''(y_{\Gamma,1}(t)) \dt\yG(t) \, \vG 
   + \iG \bigl( \fG''(y_{\Gamma,1}(t)) - \fG''(y_{\Gamma,2}(t)) \bigr) 
  \dt{y_{\Gamma,2}(t)} \, \vG
   \\
  \label{abst_var_due}
  & \quad 
  - \iG \dt\uG(t) \, \vG
    \quad \hbox{for every $(v,\vG) \in\calV$}.
\end{align}
Again, $\dt y$ possesses zero mean value. Taking into account
the previous equations at the time $s$, testing \eqref{abst_var_uno} by 
$\calN (\dt y(s))$, \eqref{abst_var_due} by $- \dt {(y,\yG)}(s)$, 
integrating over $(0,t)$ with respect to $s$, 
and adding the resulting equations leads to
\begin{align}
  \non
  & \iot \< \dt (\dt y), \calN (\dt y)>
  + \intQt \nabla (\dt w) \cdot \nabla \calN (\dt y)
  - \intQt \dt w \, \dt y
  +  \intSt \dt (\dt \yG) \dt\yG
  \\
  \non
  & \quad{} +  \intQt |\nabla \dt  y |^2
  + \intSt |\nablaG \dt\yG|^2
  = - \intQt \bigl( f''(y_1) - f''(y_2) \bigr)\dt y_2 \, \dt y
  - \intQt f''(y_1) |\dt y |^2
  \\
  & \quad{} 
  - \intSt \bigl( \fG''(y_{\Gamma,1}) - \fG''(y_{\Gamma,2}) \bigr) 
  \dt{y_{\Gamma,2}} \, \dt \yG
  - \intSt \fG''(y_{\Gamma,1}) |\dt\yG |^2
  + \intSt \dt\uG \, \dt \yG,
  \label{lhside}
\end{align}
where the first three terms can be treated, using \eqref{defN} and \eqref{simmN}, as follows
\begin{align}
  \non
  & \iot \<\dt (\dt y), \calN (\dt y)>
  + \intQt \nabla (\dt w) \nabla \calN (\dt y)
  - \intQt \dt w  \dt y
  = \frac 12 \iot \biggl( \frac d {dt} \norma{\dt y}_*^2 \biggr).
\end{align}
Integrating by parts and invoking the boundedness and the \Lip~continuity of $f''$ and $\fG''$,
we infer that
\begin{align}
  \non
  & \frac 12 \norma{\dt y(t)}^2_*
  + \frac 12 \iG |\dt \yG(t)|^2
  + \intQt |\nabla \dt y|^2
  + \intSt |\nablaG \dt \yG|^2
  \\  \non  & \quad{} 
    \leq  \frac 12 \norma{\dt y(0)}^2_*
	+ \frac 12 \iG |\dt \yG(0)|^2
  + c \intQt |y||\dt y_2||\dt y|
  + c \intQt |\dt y|^2  
  \\  & \qquad{} 
  + c \intSt |\yG||\dt y_{\Gamma,2}||\dt \yG|
  + c \intSt |\dt \yG|^2
  + c \intSt |\dt \uG||\dt \yG|,
  \label{lemma2stima}
\end{align}
where the terms on the \rhs~are denoted by $I_1,...,I_7$, in this order.
By considering the variational formulation \accorpa{primadiffdue}{secondadiff} 
and putting $t=0$, we deduce that
\begin{align}
\non
  & \<\dt y(0)  , v >
  + \iO \nabla w(0) \cdot \nabla v = 0
  \quad{} \hbox{for every $v\in V,$}
  \\ \non
  & \iO w(0) \, v
  = \iG \dt\yG(0) \, \vG
  - \iG \uG(0) \vG
  \quad{} \hbox{for every $(v,\vG) \in\calV$.}
\end{align}
Then, we test the former by $\calN (\dt y(0))$, the latter
by $-\dt (y,\yG)(0)$, and add the resulting equalities to obtain that
\begin{align}
   \non
   & \iO \dt y(0) \calN \dt y(0)
   + \iO \nabla w(0) \cdot \nabla \calN (\dt y(0)) 
   - \iO w(0) \dt y(0)
   \\
   & \qquad + \iG \dt \yG(0) \dt \yG(0)
   = \iG \uG(0) \dt \yG(0).
   \label{firstint}
\end{align}
Note that the second and third terms cancel out. Moreover, owing to 
the Young inequality we can estimate the integral on the \rhs~realizing that
\Beq
   \non
   \iO \dt y(0) \calN \dt y(0)
   + \iG |\dt \yG(0)|^2 
   = \iG \uG(0) \dt \yG(0)
   \leq
   \frac 12 \iG |\dt \yG(0)|^2
   + \frac 12 \iG |\uG(0)|^2.
\Eeq
Rearranging the terms, we deduce that
\begin{align}
   \non
   & 
   |I_1|+|I_2|
   \leq \frac 12 \iG |\uG(0)|^2
   \leq \frac 12 \norma \uG_{\C0\HG}^2
   \leq c \norma \uG_{\H1\HG}^2,
\end{align}
where the standard embedding $H^1 (0,T;\HG) \subset C^{0} ([0,T]; \HG) $
is also taken into account.
Coming back to inequality \eqref{lemma2stima}, we continue the
analysis focusing on the third integral, which can be managed as follows
\begin{align*}
   \non   
   |I_3|
   &\leq
   c \intQt |y||\dt y_2||\dt y|
   \leq c \iot \norma y _4 \norma {\dt y_2} _4 \norma {\dt y} _2
   \leq c \iot \norma y _V \norma {\dt y_2} _V \norma {\nabla \dt y} _H
   \\
   \non  
   &
   \leq \frac 14 \iot \norma {\nabla\dt y} _H^2  
   + c \norma y ^2_{\L\infty V} \norma {\dt y_2}^2_{\L2 V}
   \leq \frac 14 \iot \norma {\nabla\dt y} _H^2 
   + c \norma \uG ^2_{L^2(\Sigma)},
\end{align*}
where we applied the \andrea{\Holder, Poincar\'e and Young inequalities},
the Sobolev embedding of $V\subset \Lx4$, and at the end also the 
stability estimate \eqref{stability} along with \eqref{stab} for $y_2$. 
Moreover, combining the compactness inequality
\eqref{compactineq} with the \Poincare\ inequality \eqref{poincare} and \eqref{stability}, 
we get that
\Bsist
   |I_4|
   \leq \frac 14 \iot \norma {\nabla\dt y} _H^2 
   + c \iot \norma {\dt y} ^2_*
   \leq \frac 14 \iot \norma {\nabla\dt y} _H^2
   + c \norma \uG ^2_{L^2(\Sigma)}.
   \label{stimaB}
\Esist
The boundary terms can be dealt in a similar way as follows
\begin{align}
   \non  
    |I_5|
   &\leq c \iot \norma \yG _4 \norma {\dt y_{\Gamma,2}} _4 \norma {\dt \yG} _2
   \leq c \iot \norma \yG _{\VG} \norma {\dt y_{\Gamma,2}} _{\VG} \norma {\dt \yG} _{\HG} 
   \\
   \non
   & 
   \leq c \iot \norma {\dt \yG}_{\HG}^2 
   + c \norma {\dt y_{\Gamma,2}} ^2_{\L2\VG} \norma{ \yG}^2_{\L\infty\VG}
   \\ 
   & 
   \leq c \iot \norma {\dt \yG}_{\HG}^2 
   + c \norma \uG ^2_{L^2(\Sigma)},
   \label{stima2}
\end{align}
\Accorpa \Stime stima1 stima2   
where the fact that $y_2$ is a solution to system \Ipbl~and the inequality \eqref{stability}
turn out to be fundamental. Finally, using \eqref{stability} once more, we infer that
\Beq
   |I_7|
   \leq \norma{\dt \uG}_{L^2(\Sigma)} \norma{\dt \yG}_{L^2(\Sigma)}
   \leq c \norma{\dt \uG}_{L^2(\Sigma)} \norma \uG_{L^2(\Sigma)}
   \leq c \norma \uG_{\H1\HG}^2.
   \label{stima3}
\Eeq
Then, upon collecting the above estimates, we rearrange \eqref{lemma2stima} to realize that
\begin{align}
   \non
   & \frac 12 \norma{\dt y(t)} ^2_*
   + \frac 12 \iG |\dt \yG(t)|^2
   + \frac 12 \intQt |\nabla \dt y|^2
   + \frac 12 \intSt |\nablaG \dt \yG|^2
   \leq c \norma \uG_{\H1\HG}^2,
\end{align}
which allows us to conclude that
\Beq
   \label{partiallemma}
   \norma{(y_1,y_{\Gamma,1}) - (y_2,y_{\Gamma,2})}_{\W{1,\infty}\calG \cap \H1\calV}
   \leq c \norma{u_{\Gamma,1} - u_{\Gamma,2}}_{\H1\HG}.
\Eeq
Now, it remains to show that
\Beq
   \non
   \norma y _{\L\infty\Hdue} + \norma \yG _{\L\infty\HdueG}\leq c \norma \uG _{\H1\HG}
\Eeq 
is satisfied for some positive constant $c$.
To this aim, we test \eqref{primadiffdue} 
by $ w(t) - (w(t))^\Omega$ and integrate over $\Omega$ to get
\begin{align}
   \non
    & 
    \iO \Big|\nabla w(t)\Big|^2
   = - \< \, \dt y (t), w(t) - (w(t))^\Omega >
   \leq c \norma{\dt y (t)}_* \norma{w(t) - (w(t))^\Omega}_V
   \\
   \non
   & \quad
   \leq \frac 12 \iO \Big|\nabla w(t)\Big|^2 
   + c \norma{\dt y (t)}_*^2,
\end{align}
thanks to the \Holder, Young and \Poincare\ inequalities. Hence, 
applying \eqref{partiallemma} we find out that 
\Beq
	\non
	\norma{\nabla w}_{\L\infty H} \leq c \norma \uG _{\H1 \HG}.
\Eeq
Next, we would like to recover the full norm of $w$ in $\L\infty V$.
In this direction, we will show a bound for its mean value,
and then apply the \Poincare~inequality \eqref{poincare} to conclude.
Thus, we test equation \eqref{secondadiff} by $1$ and
integrate over $\Omega$ to obtain that
\begin{align*}
  \non
  \iO w(t)
  = \iG \dt\yG(t)
  + \iO \bigl( f'(y_1(t)) - f'(y_2(t)) \bigr)
  + \iG \bigl( \fG'(y_{\Gamma,1}(t)) - \fG'(y_{\Gamma,2}(t)) - \uG(t) \bigr),
\end{align*}
from which, owing to \eqref{poincare}, we deduce that
\Bsist
   \norma{w}_{\L\infty V}
   \leq c \norma{u_{\Gamma}}_{\H1\HG}.
   \label{stimaw}
\Esist
In order to apply a comparison principle, \andrea{we} consider
the variational formulation \accorpa{primadiffdue}{secondadiff},
and integrate by parts so to derive
the corresponding strong formulation, that holds at least 
in a distributional sense. It reads as follows
\begin{align}
   \label{strongfprima}
  &\dt y - \Delta w   = 0
   \quad \hbox{in $\, Q$},
   \\
   \label{strongseconda}
    & w = - \Delta y + f'(y_1) - f'(y_2)
   \quad \hbox{in $\, Q$}.
\end{align}
Then, comparison in \eqref{strongseconda} yields that
\Beq
	\non
  \norma{\Delta y}_{\L\infty H} 
  \leq \norma w _{\L\infty H}
  + c \norma y_{\L\infty H}
  \leq c \norma{\uG}_{\H1\HG},
\Eeq
accounting for the previous estimates, along with the regularity of $f$.
Next, \eqref{partiallemma} and the regularity theory for elliptic equation 
(see, e.g., \cite[Thms.~7.3 and~7.4, pp.\ 187-188]{LioMag}
or \cite[Thm.~3.2, p.~1.79, and Thm.~2.27, p.~1.64]{BreGil}) give us
\Beq
   \non
   \norma y _{\L\infty {H^{3/2}(\Omega)}}
   + \norma{\dn y}_{\L\infty\HG}
   \leq c \norma{\uG}_{\H1\HG},
\Eeq
which allows us to write the boundary conditions in the following form
\begin{align}
   & \dn y + \dt \yG - \DeltaG \yG + 
   \fG'(y_{\Gamma,1}) - \fG'(y_{\Gamma,2}) = \uG
   \quad \hbox{on $\, \Sigma$},
   \label{strongfultima}
   \\
   \label{strongfultimadue}
   & \dn w  = 0
   \quad \hbox{on $\, \Sigma$}.
\end{align}
Arguing in a similar manner, we can use a comparison principle in the boundary 
equation \eqref{strongfultima} to infer that
\Beq
  \non
  \norma{\DeltaG \yG}_{\L\infty\HG}
  \leq c \norma{\uG}_{\H1\HG}.
\Eeq
The boundary version of the regularity results for elliptic equations entails that 
\Beq
	\non
   \norma {\yG}_{\L\infty {H^2(\Gamma)}}
   \leq c \norma{\uG}_{\H1\HG},
\Eeq
which in turn, together with \eqref{partiallemma}, implies that
\Beq
   \norma {y}_{\L\infty {H^2(\Omega)}}
   \leq c \norma{\uG}_{\H1\HG}.
\Eeq
Then, \eqref{stability2} is completely proved.
\Edim

With these \andrea{stability results} at disposal, we are now in a position to show the
\Frechet~differentiability of $\calS$, that is, to check Theorem \ref{Fdiff}.
\andrea{Let us also point out that, owing to the different approaches employed in \cite{CGSpure},
the \Frechet\ differentiability of the control-to-state operator
was not analyzed there.}
\proof[Proof of Theorem \ref{Fdiff}]
For the sake of simplicity, we fix $\uG\in \cal U$ 
(instead of $\uopt$ used in the statement). Then, since $\cal U$ is open, 
provided we take $\hG\in \cal X$ sufficiently small, we also have that
$\uG + \hG \in \cal U$. From now on, we tacitly assume that this is the case.
Moreover, for every given $\hG\in\calX$, let us set
\Bsist
  \non
 &  \dafare{(y,\yG,w) := \hbox{solution to system \Ipbl~corresponding to $\uG$,}}
  \\ \non
  & (\yh,\yhG,\wh) := \hbox{solution to system \Ipbl~corresponding to $\uG+\hG$,}
  \\
  &\hbox{where } (y,\yG) = \calS(\uG), \ \hbox{ and } \ (\yh,\yhG) = \calS(\uG+\hG).
\Esist
For convenience, we \andrea{use} the following notation
\Bsist
  & \qh := \yh - y - \xi , \quad
  \qhG := \yhG - \yG - \xiG ,
  \aand
  \zh := \wh - w - \eta \, ,
  \non
\Esist
where $(\xi,\xiG,\eta)$ is the solution to the linearized system \Linpbl\ 
corresponding to $\hG$.
We aim to verify the \Frechet~differentiability of $\calS$ by checking the definition.
Namely, we should find a linear operator $[D{\cal S}(u_\Gamma)](\hG)$ such that
\Beq
   \non
   \calS (\uG+\hG) = \calS (\uG) + [D{\cal S}(u_\Gamma)](\hG) + o(\norma\hG_{\calX})
   \quad \hbox{in $\cal Y$} \ \  \hbox{as} \ \  \norma\hG_{\calX}\to 0\andrea{.}
\Eeq
We claim that $[D\calS(\uopt)](\hG)=(\xi,\xiG)$.
Accounting for the above notation, we realize that the above condition is equivalent 
to show that
\Beq
  \non
  \frac{\norma{(\qh,\qhG)}_{\calY}}{\norma\hG_{\calX}} \to 0 \quad \hbox{as} \quad{} 
   \norma\hG_{\calX}\to 0.
\Eeq
Furthermore, a sufficient condition consists in \andrea{proving} that
\Beq
  \norma{(\qh,\qhG)}_{\calY}
  \leq c \norma\hG_{\LS2}^2,
  \label{tesiFrechet}
\Eeq
which is the estimate we are going to \andrea{check}.
To this aim, let us consider the variational formulations for
the triplets $(\yh,\yhG,\wh)$ and $(y,\yG,w)$
satisfying problem \State\ with data $\uG+ \hG$ and $\uG$, 
and the one for $(\xi,\xiG,\eta)$ that solves the linearized system~\Linpbl.
\andrea{Then, we take the difference to obtain, \aat, that}
\begin{align}
  & \< \dt\qh(t), v >
  + \iO \nabla \zh(t) \cdot \nabla v = 0
  \quad \hbox{for every $v\in V$,}
  \label{primah}
  \\  
  & \iO \zh(t) \, v
  = \iG \dt\qhG(t) \, \vG
  + \iO \nabla\qh(t) \cdot \nabla v
  + \iG \nablaG\qhG(t) \cdot \nablaG\vG
  \non
  \\
  & \quad {}
  + \iO \bigl( f'(\yh(t)) - f'(y(t)) - f''(y(t)) \xi(t) \bigr) \, v
  \qquad
  \non
  \\
  & \quad {}
  + \iG \bigl( \fG'(\yhG(t)) - \fG'(\yG(t)) - \fG''(\yG(t)) \xiG(t) \bigr) \, \vG 
  \quad \hbox{for every $(v,\vG) \in\calV$}
  \label{secondah}
\end{align}
and that $\qh(0)=0$.
To perform our estimate, we first 
add to both sides of \eqref{secondah} the term $\iO \qh(t) \, v $ and 
the corresponding boundary contribution $\iG \qhG(t) \, \vG $.
Then, we test \eqref{primah} and this new \eqref{secondah}, written at the time $s$,
by $\calN(\dt\qh(s))$ and $-\dt(\qh,\qhG)(s)$, respectively.
Adding the resulting equalities and integrating over $(0,t)$ for an arbitrary $t\in(0,T)$,
leads~to
\begin{align}
  & \iot \<\dt\qh , \calN(\dt\qh)> 
  + \intQt \nabla\zh \cdot \nabla\calN(\dt\qh)
  - \iot \<\dt\qh , \zh >
  + \intSt |\dt\qhG|^2
  \non
  \\
  & \quad\quad  
  \non
    + \frac 12 \normaV {\qh(t)}^2
  + \frac 12 \norma {\qhG(t)}^2_{\VG}= - \iot \< \dt\qh , f'(\yh) - f'(y) - f''(y) \xi - \qh>
  \\
  & \quad\quad
  - \intSt \bigl( \fG'(\yhG) - \fG'(\yG) - \fG''(\yG) \xiG - \qhG \bigr) \, \dt\qhG.
  \label{perFrechet}
\end{align}
As before, the first three integrals on the \lhs\ can be easily handled
with a cancellation, so that
\Bsist
  \non
  \iot \< \dt\qh , \calN(\dt\qh)>
  + \intQt \nabla\zh \cdot \nabla\calN(\dt\qh)
  - \iot \<\dt\qh , \zh >
  = \iot \norma{\dt \qh}_*^2
  \geq 0.
\Esist
Note that the other terms of the \lhs~are nonnegative.
Owing to the regularity of the potentials, we can invoke
the Taylor \andrea{formula} with integral remainder \andrea{for the function $f'$}
at $y$. Recalling that
$\yh-y=\xi+\qh$, we have
\Beq
  f'(\yh) - f'(y) - f''(y) \xi
  = f''(y) \, \qh 
  +  \int_0^1 f'''(y + \zeta (\yh -y))(1- \zeta) \andrea{(\yh -y)^2} d\zeta.
  \label{taylorinteg}
\Eeq
As the \rhs~of \eqref{perFrechet} is concerned, let us estimate the first integral
as follows
\begin{align}
   \non
   & - \iot \< \dt \qh ,  f'(\yh) - f'(y) - f''(y) \xi - \qh >
   \\ \non
   & \quad
   \leq \frac 12 \iot \norma{\dt \qh}_*^2
   + \,c \iot \Big\|{ f'(\yh) - f'(y) - f''(y) \xi - \qh }\Big\|^2_V 
   \\ 
   &\quad
   \leq \frac 12 \iot \norma{\dt \qh}_*^2   
   + \, c \iot \normaV{\qh }^2
   \non\\ 
   &\qquad
     + \,c \iot  \Big\|{f''(y) \, \qh  
   +\int_0^1 f'''(y + \zeta (\yh -y))(1- \zeta)(\yh -y)^2 d\zeta
   } \, \Big\|_V^2,
\end{align}
thanks to \eqref{young} and \eqref{taylorinteg}.
Moreover, the last term can be dealt as follows
\begin{align}
   \non
   & c\iot \Big\|{f''(y) \, \qh  +  \int_0^1 f'''(y + \zeta (\yh -y))(1- \zeta)(\yh -y)^2 d\zeta}\Big\|_V^2
   \\ \non
   & \quad{} 
      \leq c\iot \normaV{f''(y) \, \qh}^2 
	+ c\iot \Big\|{ \int_0^1 f'''(y + \zeta (\yh -y))(1- \zeta)(\yh -y)^2 d\zeta}\Big\|_V^2 .
\end{align}
Proceeding with a separate analysis, we obtain that
\begin{align}
   \non
    \iot \Big\|{f''(y) \, \qh}\Big\|_V^2 
   &=  \iot \biggl( \normaH{f''(y) \, \qh}^2 + \normaH{\nabla(f''(y) \, \qh)}^2 \biggr)
   \\
   \non
   & 
   \leq \iot \biggl( 
   \norma{f''(y)}_\infty^2 \normaH{\qh}^2 
   + \norma{f''(y)}_\infty^2 \normaH{\nabla\qh}^2 
   + \norma{f'''(y)}_\infty^2 \normaH{\nabla y \cdot \qh}^2 
   \biggr)
   \\
   \non
   & 
    \leq c \iot \normaV{\qh}^2
   + c \iot \norma{\nabla y}_4^2\norma{\qh}_4^2
    \leq c \iot \normaV{\qh}^2
   + c \iot \norma{\nabla y}_V^2\norma{\qh}_V^2
   	\\ \non
	&  
   \leq 
   c \iot \normaV{\qh}^2
   + c \norma y _{\L\infty \Hdue}^2 \iot \normaV{\qh}^2 ,
\end{align}
owing to the fact that $y$, as a solution to \Pbl, satisfies \eqref{stab}
and thanks to the Sobolev embedding $V\subset \Lx4$. 
\dafare{
Furthermore, the second part can be handled as follows
\begin{align*}
	&
	\iot \Big\|{ \int_0^1 f'''(y + \zeta (\yh -y))(1- \zeta)(\yh -y)^2 d\zeta}\Big\|_V^2 
	\\
	& \quad
	\leq
	\iot \int_0^1 \Big\|{  f'''(y + \zeta (\yh -y))(1- \zeta)(\yh -y)^2 }\Big\|_V^2 d\zeta	
	\\
	& \quad
	\leq
	\iot \int_0^1 \Big\|{  f'''(y + \zeta (\yh -y))(1- \zeta)(\yh -y)^2 }\Big\|_H^2 d\zeta
	\\ & \qquad\qquad
	+\iot \int_0^1 \Big\|{\nabla\Bigl(  
				f'''(y + \zeta (\yh -y))(1- \zeta)(\yh -y)^2 \Bigr)}\Big\|_H^2 d\zeta
	\\
	& \quad
	\leq 
	\sup_{0\leq\zeta\leq 1} \big\| f'''(y + \zeta (\yh -y)) \big\|_\infty
	\iot \Big\|(\yh -y)^2 \Big\|_H^2
	\\ & \qquad\qquad
	+
	\iot \int_0^1 \Big\|{  f^{(iv)}(y + \zeta (\yh -y))\nabla (y + \zeta (\yh -y))(1- \zeta)
	(\yh -y)^2 }\Big\|_H^2 d\zeta	
	\\ & \qquad\qquad
	+
	\iot \int_0^1 \Big\|{  f'''(y + \zeta (\yh -y))(1-\zeta)2(\yh -y)\nabla (\yh -y)}\Big\|_H^2 d\zeta.
\end{align*}
Consequently, we have that 
\begin{align*}
	& 	\iot \Big\|{ \int_0^1 f'''(y + \zeta (\yh -y))(1- \zeta)(\yh -y)^2 d\zeta}\Big\|_V^2 
	\\
	& \quad
	\leq
	c \iot \Big\|(\yh -y)^2 \Big\|_H^2
	+\sup_{0\leq\zeta\leq 1} \big\| f^{(iv)}(y + \zeta (\yh -y)) \big\|_\infty
	\iot \Big\|(|\nabla y| +|\nabla \yh|) (\yh -y)^2\Big\|_H^2
	\\ & \qquad\qquad
	+2 \sup_{0\leq\zeta\leq 1} \big\| f'''(y + \zeta (\yh -y)) \big\|_\infty
	\iot  \Big\| (\yh -y)\nabla (\yh -y)\Big\|_H^2
	\\ & \quad
	\leq
	c \iot \Big\|\yh -y \Big\|_{4}^4
	+ c \iot \Bigl(\norma{\nabla y}_6^2 +\norma{\nabla \yh}_6^2\Bigr) \Big\|\yh -y\Big\|_6^4
	+ c \iot \Big\|\yh -y\Big\|_4^2 \Big\|\nabla (\yh -y) \Big\|_4^2
	\\ & \quad
	\leq
	c \iot \Big\|\yh -y \Big\|_{V}^4
	+ c \iot \Bigl(\norma{\nabla y}_V^2 +\norma{\nabla \yh}_V^2\Bigr) \Big\|\yh -y\Big\|_V^4
	\\ & \qquad\qquad
	+ c \iot \Big\|\yh -y\Big\|_V^2 \Big\|\nabla (\yh -y) \Big\|_V^2
	\\ & \quad
	\leq{}
	 c \norma{\hG}^4_{\H1\HG},
\end{align*}
where the Sobolev embeddings $V \subset \Lx4$ and $V \subset \Lx6$, 
\andrea{and}
the stability estimate \eqref{stability2} have been used along with
the fact that $y$ and $\yh$, as solutions
to \Ipbl, satisfy \eqref{stab}.}
Summarizing, we have just shown that
\Beq
	\non
	c \iot \Big\|{f''(y) \, \qh  + \int_{y}^{\yh} f'''(\gamma)(\yh-\gamma)^2 d\gamma} \Big\|_V^2
	\leq c \iot \normaV{\qh}^2
	+ c \norma{\hG}^4_{\H1\HG}.
\Eeq
\dafare{
Using the Taylor formula corresponding to \eqref{taylorinteg} 
for the the nonlinearity $\fG'$,
combined with the Young inequality and the stability estimate \eqref{stability},
we control the last term of \eqref{perFrechet} by
\begin{align}
   \non
   & - \intSt \bigl( \fG'(\yhG) - \fG'(\yG) - \fG''(\yG) \xiG - \qhG \bigr) \, \dt\qhG
    \\
   \non
   & \quad{}
   \leq \frac 12 \intSt |\dt \qhG|^2
   + \frac 12 \intSt | \fG'(\yhG) - \fG'(\yG) - \fG''(\yG) \xiG - \qhG|^2
   \\
   \non
   & \quad{}
   \leq \frac 12 \intSt |\dt \qhG|^2
      \\ \non & \qquad \qquad
       {} + \frac 12 \intSt \Big| 
        \fG''(\yhG)\qhG
	+  \int_0^1   \fG '''\bigl(\yG + \zeta (\yhG -\yG)\bigr)(1- \zeta) \andrea{(\yhG -\yG)^2} d\zeta
   - \qhG
   \Big|^2
   \\
   \non
   & \quad{}
   \leq \frac 12 \intSt |\dt \qhG|^2
   + c \intSt |\qhG|^2
   + c \norma{\fG''(\yhG)}_\infty^2\intSt | \qhG|^2
   \\ \non & \qquad \qquad
   + c \sup_{0\leq\zeta\leq 1} \Big\|\fG '''\bigl(\yG + \zeta (\yhG -\yG)\bigr) \Big\|_\infty
   \intSt |\yhG -\yG|^4
   \\
   \non
   & \quad{}
   \leq \frac 12 \intSt |\dt \qhG|^2
   + c \intSt |\qhG|^2
   + c \norma{\hG}_{L^2(\Sigma)}^4.
\end{align}}
Summing up, upon collecting all the above estimates, we realize that the inequality
\begin{align}
   \non
   & \frac 12 \iot \norma{\dt \qh}_*^2 
   + \frac 12 \intSt |\dt \qhG|^2
   + \frac 12 \normaV{\qh(t)}^2
   + \frac 12 \norma{\qhG(t)}^2_{\VG}
   \\
   \non
   & \quad{} 
   \leq c \iot \normaV{\qh}^2
   + c \intSt |\qhG|^2
   + c \norma{\hG}^{4}_{\H1\HG},
\end{align}
has been proved
\andrea{whence}
a Gronwall argument directly yields \eqref{tesiFrechet}. 
\qed

\section{Optimality conditions}
\label{OPTIMALITY}
\setcounter{equation}{0}
\subsection{The adjoint system}

This section is completely devoted to the investigation of the adjoint
system and to the necessary conditions for optimality. 
Let us begin with the task of ensuring the well-posedness of system \Pbladj, that
is checking Theorem \ref{adjexist}.
\andrea{Before \pier{going on}, it is worth mentioning here that
the technique below differs from the one employed in \cite{CGS_OPT}.
Indeed, the key argument to prove the well-posedness of the adjoint problem
in \cite{CGS_OPT} relies on \pier{the interpretion of} the adjoint system
as a suitable abstract Cauchy problem in a general mathematical framework\pier{:
then, after proving that the involved operators verify
some properties like coercivity and continuity,
the existence and uniqueness of the solution are deduced from some classical results.}
In this direction, let us emphasize that the \pier{outlined analysis suggested the 
authors to confine themselves} to 
the investigation to the case $b_\Omega=b_\Gamma=0$
(see also \cite[Rem.~5.6]{CGS_OPT}, where a 
possible way to overcome this restriction is explained 
\pier{by involving} weighted Lebesgue spaces).}

\proof[Proof of Theorem \ref{adjexist}]
We will tackle the proof in two steps. 
In the first one, we will check the existence
of a solution with the required regularity, whereas in the second step, we
will point out that such a solution is indeed unique.
From now on, let us convey that
$\uopt$ and $(\yopt,\yGopt)$ stand for an optimal control with the 
corresponding optimal state, respectively. 

\noindent
{\bf Existence}
The key idea for the existing part is showing that system \Pbladj\ 
can be rewritten as an initial boundary value problem
which complies with the framework of \cite[Thm.~2.3, p.~977]{CGS}.

Moreover, since we are going to reverse the time 
with the following change of variable $t\mapsto T-t$, it turns
out to be useful to set
\begin{align*}
	\non
	\widetilde{\phi}_Q(t)&:=\phQ(T-t), \ \ 
	\widetilde{\phi}_\Sigma(t):=\phS(T-t), 
	\widetilde{q}(t):=q(T-t), \ \ \widetilde{p}(t):=p(T-t), \ \ 
	\\
	\widetilde{\lam}(t)&:= \lam(T-t), \ \ 
	\widetilde{\lam}_\Gamma(t):={\lam}_\Gamma(t), \ \
	\widetilde{q}_{\Gamma}(t):=\widetilde{q}\suG (t) \ \ \aat.
\end{align*}
\andrea{Therefore, after substituting} $t$ with $T-t$\andrea{, we} realize that
system \Pbladj\ can be reformulated
as the initial boundary value problem 
\begin{align}
	& \label{Agg_back_prima}
	\iO \widetilde{q}(t) v
	=
	\iO \nabla \widetilde{p}(t) \cdot \nabla v 
	\quad \hbox{for every $v \in V, \ \aat$},
	\\
	\non
	&
	\iO \dt \widetilde{p} (t) v
	+ 	\iO \nabla \widetilde{q}(t) \cdot \nabla v 
	+ 	\iO \widetilde{\lam}(t)\widetilde{q}(t) v 
	+
	\iG \dt \widetilde{q}_\Gamma (t) v_\Gamma
	+ 	\iG \nabla_\Gamma \widetilde{q}_\Gamma(t) \cdot \nabla_\Gamma v _\Gamma
	\\  
	\non
	& \quad\quad
	+ 	\iG \widetilde{\lam}_\Gamma(t)\widetilde{q}_\Gamma(t) v _\Gamma
	= 
	\iO \widetilde{\phi}_Q(t) v
	+ \iG \widetilde{\phi}_\Sigma(t) v_\Gamma
	\\
	& \hspace{7cm}
	 \hbox{for every $(v,\vG) \in \cal V, \ \aat$},
	\label{Agg_back_seconda}
	\\ &
	\iO \widetilde{p}(0) v
	+ 	\iG \widetilde{q}_\Gamma(0) v_\Gamma
	=
	\iO {\phi_\Omega} v
	+ \iG {\phi_\Gamma} v_\Gamma
	\quad \hbox{for every $(v,\vG) \in \cal V$}.
	\label{Agg_back_terza}
\end{align}
\Accorpa\Aggback Agg_back_prima Agg_back_terza
We claim that \Aggback~can be studied with the help of
\cite[Thm.~2.3]{CGS}. 
In this direction, let us proceed indirectly. Hence, we pick a
function $\Phi\in V$ such that $(\Phi, \phi_\Gamma)\in \cal V$, i.e.
$\Phi \suG= \phi_\Gamma$.
Then, we take into account the problem of looking for a triplet 
$(r,r_\Gamma,\mu)$ which satisfies the following problem:
\begin{align}
	& \label{pier2}
	 r_\Gamma(t)= r(t)\suG \quad \aat,
	\\
	 \label{Agg_ausiliario_prima}
	\< \dt r(t), v>
	&=
	\iO \nabla \mu(t) \cdot \nabla v 
	\quad \hbox{for every $v \in V, \ \aat$},
	\\
	\non	
	\iO \mu (t) v
	&=
 	\iO \nabla r(t) \cdot \nabla v 
	+ 	\iO \widetilde{\lam}(t)r(t) v 
	+
	\iG \dt r_\Gamma (t) v_\Gamma
	+ 	\iG \nabla_\Gamma r_\Gamma(t) \cdot \nabla_\Gamma v _\Gamma
	\\ 
	\non
	& \quad
	+ 	\iG \widetilde{\lam}_\Gamma(t)r_\Gamma(t) v _\Gamma
	= 
	\iO \widetilde{\phi}_Q(t) v
	+ \iG \widetilde{\phi}_\Sigma(t) v_\Gamma
	\\ &\hspace{5cm}
	\quad \hbox{for every $(v,\vG) \in \cal V, \ \aat$},
	\label{Agg_ausiliario_seconda}
	\\
	r(0)&= \Phi
	\quad \hbox{in $\Omega$}, 
	\label{Agg_ausiliario_terza}
\end{align}
\Accorpa\Aggausiliario pier2 Agg_ausiliario_terza
where the functions 
$\widetilde{\lam},\widetilde{\lam}_\Gamma,\widetilde{\phi}_Q,\widetilde{\phi}_\Sigma$
are the same as above. Furthermore, the previous investigation,
along with \eqref{hpzzzz}, leads us to realize that
\begin{align*}
	 \widetilde{\lam}&\in L^\infty(Q) \cap \L\infty {W^{1,3}(\Omega)},
 	\ \
	\widetilde{\lam}_\Gamma\in L^\infty(\Sigma) ,
	\\
	\widetilde{\phi}_Q &\in \H1 H,
 	\ \
	\widetilde{\phi}_\Sigma \in \L2 {\HG},
	\ \
	(\Phi, \phi_\Gamma) \in \cal V.
\end{align*}
Therefore, the assumptions of \cite[Thm.~2.3, p.~977]{CGS} are satisfied so that
the existence of a triplet $(r,r_\Gamma,\mu)$\pier{,}
which solves \Aggausiliario\ and enjoys the following regularity
\begin{align*}
   (r,r_\Gamma ) &\in \H1 {{\cal G}} \cap \L\infty {{\cal V}}  \cap \L2 {\Hx2 \times H^2(\Gamma)},
  \\
    \mu &\in \L2 V,
\end{align*}
directly follows.
We are then reduced to show that system \Aggback\ can be written
in the form of \Aggausiliario. We claim that the following choice realizes
this goal:
\begin{align}
	\label{deftilde}
	\widetilde{q}:= r,
	\ \ 
	\widetilde{q}_\Gamma:= r_\Gamma, 
	\ \ 
	\widetilde{p}(t):= \phi_\Omega - \iot \mu(s) ds \quad \aat.
\end{align}
In fact, by differentiating the last term, we deduce that 
$\mu = - \dt \widetilde{p} \ \ \aeQ$, so that
\eqref{Agg_ausiliario_seconda} implies \eqref{Agg_back_seconda}. Moreover,
integrating \eqref{Agg_ausiliario_prima} with respect to $t$ and using
\eqref{Agg_ausiliario_terza} yield
\begin{align*}
	\iO r(t) v 
	+ \iO \nabla \iot \mu(s) ds \cdot \nabla v
	= 
	\iO \Phi \,  v
	\quad \hbox{for every $v \in V,$}
\end{align*}
which, owing to \eqref{deftilde}, entails that
\begin{align*}
	\iO \widetilde{q}(t) v 
	+ \iO \nabla (- \widetilde{p}(t) 
	+ \phi_\Omega) 
	\cdot \nabla v
	= 
	\iO \Phi \,  v
	\quad \hbox{for every $v \in V.$}
\end{align*}
Hence, provided we require that 
\begin{align}
	\label{condition_prima}
	\iO \nabla \phi_\Omega \cdot \nabla v
	= 
	\iO \Phi \,  v
	\quad \hbox{for every $v \in V,$}
\end{align}
\eqref{Agg_back_prima} follows from \eqref{Agg_ausiliario_prima}.
Besides, \eqref{pier2}, \eqref{Agg_ausiliario_terza} and \eqref{deftilde} imply that
\begin{align*}
	\label{}
	\widetilde{p}(0)=\phi_\Omega,
	\ \ 
	\widetilde{q_\Gamma}(0)= \phi_\Gamma \quad \hbox{in $\Omega$},
\end{align*}
whence \eqref{Agg_back_terza} immediately follows by testing by $(v,\vG)\in\cal V$ and
integrating over $\Omega$.
Summing up, equation \eqref{condition_prima} gives, in turn, that $(\Phi)^\Omega=0 $ and 
also that $\phi_\Omega$ solves the following elliptic problem
$$\begin{cases}
	- \Delta \phi_\Omega &= \Phi \quad \hbox{in $\Omega$},
	\\
	\dn \phi_\Omega &= 0 \quad \hbox{on $\Gamma,$}
\end{cases}$$
which entails that $\phi_\Omega={\cal N }(\Phi) + (\phi_\Omega)^\Omega$.
If all these compatibility conditions on $\Phi$ and $\phi_\Omega$ are in force, 
we have just checked that system \Aggback\ can be rewritten in the form of \Aggausiliario.
Thus, owing to \cite[Thm.~2.3, p.~977]{CGS}, there exists a triplet
$(\widetilde{q},\widetilde{q}_\Gamma,\widetilde{p})$, which
solves \Aggback\ and possesses the following regularity
\begin{align*}
   (\widetilde{q},\widetilde{q}_\Gamma) &\in \H1 \GO 
  	\cap \L\infty \VO \cap \L2 {\Hx2 \times H^2(\Gamma)},
  \\
  \widetilde{p}  & \in \H1 V.
\end{align*}
\dafare{
Lastly, owing to the above regularity,
along with comparison in the strong formulation of \eqref{primaadjpure},
we easily realize that $\Delta \widetilde p \in \L2{\Hx2}$, so that
the elliptic regularity theory ensures that $ \widetilde p \in \L2{\Hx4}$.
}
\Brem
Let us point out that in \cite{CGS_OPT}, where the analogous control problem 
for the viscous case was treated, the conditions $\bO=\bG=0$ have been
required in order to handle the adjoint system.
Note that this restriction leads to consider
$\phO = 0$ in $\Omega$, $\phG=0$ in $\Gamma$, $\Phi=0$ in $\Omega$, which surely 
fulfill our requirements.
\Erem

\noindent
{\bf Uniqueness}
We proceed by contradiction assuming the existence of, at least,
two solutions $(\widetilde{q}_i,\widetilde{q}_{\Gamma\!,i},\widetilde{p}_i)$,
$i=1,2$, to system \Pbladj. Then, we set 
\Beq
	\non
     \widetilde{q} := \widetilde{q}_1 - \widetilde{q}_2, \quad
     \widetilde{q}_\Gamma := \widetilde{q}_{\Gamma\!,1} - \widetilde{q}_{\Gamma\!,2}, \quad
	 \widetilde{p} := \widetilde{p}_1 - \widetilde{p}_2,
\Eeq
and we are going to show that the only possibility is  
$\widetilde{p}=\widetilde{q}=\widetilde{q}_\Gamma=0$.
In this direction, we write system \Aggback\ for both the solutions 
$(\widetilde{q}_i,\widetilde{q}_{\Gamma\!,i},\widetilde{p}_i)$,
$i=1,2$, and take the difference. 
Note that, taking $(v,0)\in \cal V$ in \eqref{Agg_back_terza},
we get $\widetilde{p}_1(0)=\widetilde{p}_2(0)=\phO$ in $\Omega$ and 
by comparison also that $\widetilde{q}_{\Gamma\!,1}(0)=\widetilde{q}_{\Gamma\!,2}(0)=\phG$.
Thus, we have that
\begin{align}
	& \label{Agg_diff_prima}
	\iO \widetilde{q}(t) v
	=
	\iO \nabla \widetilde{p}(t) \cdot \nabla v 
	\quad \hbox{for every $v \in V,\ \aat$},
	\\
	\non
	&
	\iO \dt \widetilde{p} (t) v
	+ 	\iO \nabla \widetilde{q}(t) \cdot \nabla v 
	+ 	\iO \widetilde{\lam}(t)\widetilde{q}(t) v 
	+
	\iG \dt \widetilde{q}_\Gamma (t) v_\Gamma
	+ 	\iG \nabla_\Gamma \widetilde{q}_\Gamma(t) \cdot \nabla_\Gamma v _\Gamma
	\\ 
	& \quad
	+ 	\iG \widetilde{\lam}_\Gamma(t)\widetilde{q}_\Gamma(t) v _\Gamma
	= 
	0
	\quad \hbox{for every $(v,\vG) \in \cal V,\ \aat$},
	\label{Agg_diff_seconda}
	\\ &
	\widetilde{p}(0)= 0, \ \ \widetilde{q}_\Gamma(0)= 0 \quad \hbox{in $\Omega$}.
	\label{Agg_diff_terza}
\end{align}
Next, we test equation \eqref{Agg_diff_prima} 
by $-\dt \widetilde{p}$, \eqref{Agg_diff_seconda} by 
$(\widetilde{q},\widetilde{q}_{\Gamma})$,
and \eqref{Agg_diff_prima} once more by $K \widetilde{q}$, for a constant $K$, yet
to be determined. 
Summing the obtained equalities and rearranging the terms lead to
\begin{align*}
	\non
	&\frac 12 \frac d {dt} \iO |\nabla \widetilde{p}|^2
	+ K \iO  |\widetilde{q}|^2
	+ \iO |\nabla \widetilde{q}|^2
	+	\frac 12 \frac d {dt} \iG |\widetilde{\qG}|^2
	+ \iG |\nablaG \widetilde{\qG}|^2
	\\ & \quad
	=
	K \iO \nabla \widetilde{p} \cdot \nabla \widetilde{q}
	- \iO \widetilde{\lam} |\widetilde{q}|^2 
	- \iG \widetilde{\lam}_\Gamma|\widetilde{\qG}|^2 
	\quad \hbox{$a.e.$ in $(0,T)$,}
\end{align*}
where the integrals on the \rhs~are denoted by $I_1,I_2$ and $I_3$,
respectively. 
Using the Young inequality
and the boundedness of $\widetilde{\lam}_\Gamma$, we deduce that
\begin{align*}
	|I_1| + |I_3|
	\leq 
	\frac 12 \iO |\nabla \widetilde{q}|^2
	+\frac {K^2} 2 \iO |\nabla \widetilde{p}|^2
	+ c \iG |\widetilde{\qG}|^2
	\quad \hbox{$a.e.$ in $(0,T)$.}
\end{align*}
Moreover, the boundedness of $\widetilde{\lam}$ allows us to infer that
\begin{align*}
	|I_2| \leq 
	\norma{\widetilde{\lam}}_{\infty} \iO  |\widetilde{q}|^2,
\end{align*}
and we move it to the \lhs.
Finally, we rearrange the terms and integrate over $(0,t)$ to obtain that
\begin{align*}
	\non
	&\frac 12 \iO |\nabla \widetilde{p}(t)|^2
	+ (K - \norma{\widetilde{\lam}}_{\infty}) \intQt  |\widetilde{q}|^2
	+ \frac 12 \intQt |\nabla \widetilde{q}|^2
	+	\frac 12 \iG |\widetilde{\qG}(t)|^2
	+ \intSt |\nablaG \widetilde{\qG}|^2
	\\ & \quad
	\leq
	\frac {K^2} 2 \intQt |\nabla \widetilde{p}|^2
	+c \intSt |\widetilde{q}_\Gamma|^2
	\quad \hbox{for all $t \in (0,T)$.}
\end{align*}
Hence, taking the constant $K$ \andrea{large} enough, we apply
the Gronwall lemma to conclude that
\begin{align*}
	\norma{\nabla \widetilde{p}}_{\L\infty H}
	+
	\norma{\widetilde{q}}_{\L2 V}
	+
	\norma{\widetilde{q}_\Gamma}_{\L\infty H\cap \L2 V}
	\leq 0,
\end{align*}
which yields
\begin{align*}
	\widetilde{q}=0,\ \  \widetilde{q}_{\Gamma}=0, \ \ \nabla\widetilde{p}=0.
\end{align*}
Hence, we realize that 
$\widetilde{p}$ has to be constant with respect to the space variable.
On the other hand, comparison in \eqref{Agg_diff_seconda} produces $\dt \widetilde{p}=0$, 
so that $\widetilde{p}$ has also to be constant in time
\andrea{and such that} \eqref{Agg_diff_terza} \andrea{is verified.
Therefore, we infer that} $\widetilde p = 0$ and
Theorem \ref{adjexist} is completely proved.
\qed

\subsection{Necessary optimality conditions}

The final step of \andrea{the} work consists in proving Theorem~\ref{CNoptadj}
by deriving the first-order optimality conditions. This will also point out that
\Pbladj\ yields the adjoint system for \Pbl.

\Bprop
\label{CNopt}
Let $\uopt$ and $(\yopt,\yGopt)$ be an optimal control with the corresponding 
state. Then, inequality \eqref{cnopt} holds true.
\Eprop
\Bdim
In order to prove \eqref{cnopt}, we essentially make use of \eqref{precnopt}. 
In fact, we make explicit \eqref{precnopt}
exploiting the \Frechet\ differentiability of $\calS$ and 
the chain rule. As a matter of fact, denoting by
$\widetilde\calS:\calU\to\calY\times\calX$ the function defined by
$\widetilde\calS(\uG):=(\calS(\uG),\uG)$, we realize that
Theorem~\ref{Fdiff} yields
\Beq
  D\widetilde\calS(\uG):
  \hG \mapsto \bigl( [D\calS(\uG)](\hG),\hG \bigr)
  = (\xi,\xiG,\hG)
  \quad \hbox{for the admissible $\hG\in\calX$},
  \non
\Eeq
where $(\xi,\xiG,\eta)$ is the solution to the linearized system 
\Linpbl\ corresponding to~$\hG$.
On the other hand, if we consider the cost functional $\calJ $
as a mapping from $\calY\times\calX$ to~$\erre$,
its \Frechet\ derivative at $(y,\yG,\uG)\in\calY\times\calX$
is \sfw ly given by
\begin{align*}
  & [D\calJ(y,\yG,\uG)](k,k_\Gamma,\hG)
  = \bQ \intQ (y-\zQ) k 
  + \bS \intS (\yG-\zS) k_\Gamma
  + \bO \iO (y(T)-\zO) k(T) \non
  \\
  & \quad
  + \bG \iG (y(T)-\zG) k_\Gamma(T)
  + \bz \intS \uG \hG 
  \quad \hbox{for $(k,k_\Gamma)\in\calY$ and $\hG\in\calX$}.
  \non
\end{align*}
Hence, since $\redJ=\calJ\circ\widetilde\calS$, 
the chain rule implies that 
\begin{align*}
	 [D\redJ(\uG)](\hG)
  &= [D\calJ(\widetilde\calS(\uG)] \bigl( [D\widetilde\calS(\uG)] (\hG) \bigr)
  = [D\calJ(y,\yG,\uG)] (\xi,\xiG,\hG)
  \non
  \\
  &=  \bQ \intQ (y-\zQ) \xi 
  + \bS \intS (\yG-\zS) \xiG
  + \bO \iO (y(T)-\zO) \xi(T)
  \non
  \\
  & \quad\quad {}
  + \bG \iG (y(T)-\zG) \xiG(T)
  + \bz \intS \uG \hG.
  \non
\end{align*}
Therefore, \eqref{cnopt} immediately follows from~\eqref{precnopt} by choosing in the
above calculations $(y,\yG,\uG)=(\yopt,\yGopt,\uopt)$.
\Edim

\proof[Proof of Theorem \ref{CNoptadj}]
For the sake of simplicity, we will avoid writing explicitly the time variable
in the calculations below.
Moreover, for the reader's convenience, we rewrite the variational formulation 
of the linearized and the adjoint system, respectively. They read as follows
\begin{align*}
  & 
  - \iO \dt\xi  v
  - \iO \nabla\eta \cdot \nabla v = 0
  \quad \hbox{for every $v\in V$},
  \\
    & \iO \eta v
  = \iG \dt\xiG \, \vG
  + \iO \nabla\xi \cdot \nabla v
  + \iG \nablaG\xiG \cdot \nablaG\vG
  + \iO \lam \, \xi \, v
  + \iG \bigl( \lamG \, \xiG - \hG \bigr) \, \vG
  \non
  \\
  & \hspace{10.7cm} \hbox{for every $(v,\vG) \in\calV$},
  \\ & 
  -\iO q \, v
  + \iO \nabla p \cdot \nabla v
  =0 
  \quad \hbox{for every $v\in V,$}
  \\
  & - \iO \dt p \, v
  + \iO \nabla q \cdot \nabla v
  + \iO  \lam q  v
  - \iG \dt \qG \, \vG
  + \iG \nablaG  \qG \cdot \nablaG\vG
  + \iG  \lamG \qG \vG
  \non
  \\
  &  
  \hspace{3.8 cm}
  = 
  \bQ \iO (\yopt - \zQ) v
  + \bS \iG (\yGopt - \zS) \vG
  \quad \hbox{for every $(v,\vG) \in\calV$},
\end{align*}
with the corresponding initial conditions
\begin{align*}
	(\xi,\xiG)(0)=(0,0) \quad \hbox{in $\Omega$},
\end{align*}
and final conditions
\begin{align*}
   \non \iO p(T) v 
   + \iG \qG(T) \vG
   = \bO \iO (\yopt(T) - \zO) v(T) 
   + \bG \iG (\yGopt(T) - \zG) \vG(T)
   \\
   \quad \hbox{for every $(v,\vG) \in\calV$},
\end{align*}
respectively.
Then, we test these formulations 
by $p, (q, \qG), \eta$ and $(\xi,\xi_\Gamma)$, in this order.
Adding the resulting equalities, integrating over $(0,t)$
and by parts, and using the initial condition for $\xi$ 
and the final ones for $p$ and $\qG$, lead us to infer that the most of the 
terms cancel out and it remains
\begin{align*}
	& \intS \qG \hG 
  = \intQ \bQ (\yopt - \zQ) \xi
  + \intS \bS (\yopt_\Gamma - \zS) \xiG
  \non
  \\
  & \quad{}
  + \iO \bO (\yopt(T) - \zO) \xi(T)
  + \iG \bG (\yGopt(T) - \zG) \xiG(T),
\end{align*}
which is the desired conclusion since it allows
us to obtain \eqref{cnoptadj}, where $\hG = \vG -\uopt$, from \eqref{cnopt}.
\qed

\subsection*{Acknowledgments}
\dafare{The current contribution originated from the work done by
Andrea Signori for the preparation of his master thesis discussed
at the University of Pavia on September 2017.
Actuallly, the paper turns out to offer some improvement on the 
results there contained.
The research of Pierluigi Colli is supported by the Italian Ministry of Education,
University and Research~(MIUR): Dipartimenti di Eccellenza Program (2018--2022)
-- Dept.~of Mathematics ``F.~Casorati'', University of Pavia.
In addition, PC gratefully acknowledges some other
support from the GNAMPA (Gruppo Nazionale per l'Analisi Matematica,
la Probabilit\`a e le loro Applicazioni) of INdAM (Istituto
Nazionale di Alta Matematica) and the IMATI -- C.N.R. Pavia, Italy. }

\footnotesize
\vspace{3truemm}

\Begin{thebibliography}{10}


\bibitem{BreGil}
F. Brezzi and G. Gilardi,
Part~1:
Chapt.~2, Functional spaces, 
Chapt.~3, Partial differential equations,
in ``Finite element handbook'',
H. Kardestuncer and D. {H.\ Norrie}  eds.,
McGraw-Hill Book Company, NewYork, 1987.

\bibitem{CH}
J. W. Cahn and J. E. Hilliard, 
Free energy of a nonuniform system. I. Interracial free energy, 
{\it J. Chem. Phys.} {\bf 28} (1958), 258-267. 

\bibitem{CaCo}
L. Calatroni and P. Colli,
Global solution to the Allen-Cahn equation with singular potentials and dynamic boundary conditions,
{\it Nonlinear Anal.} 
{\bf 79} (2013), 12-27.

\bibitem{CFP} 
R. Chill, E. Fa\v sangov\'a and J. Pr\"uss,
Convergence to steady states of solutions of the \CH~equation with dynamic boundary conditions,
{\it Math. Nachr.} 
{\bf 279} (2006), 1448-1462.

\pier{\bibitem{CFGS1}
P. Colli, M. H. Farshbaf-Shaker, G. Gilardi and J. Sprekels,
Optimal boundary control of a viscous Cahn--Hilliard system with dynamic boundary condition and double obstacle potentials, 
{\it SIAM J. Control Optim.} {\bf 53} (2015), 2696-2721.
\bibitem{CFGS2}
P. Colli, M. H. Farshbaf-Shaker, G. Gilardi and J. Sprekels,
Second-order analysis of a boundary
control problem for the viscous Cahn--Hilliard equation with dynamic boundary conditions,
{\it Ann. Acad. Rom. Sci. Ser. Math. Appl.} {\bf 7} (2015), 41-66.
\bibitem{CFS}
P. Colli, M. H. Farshbaf-Shaker and J. Sprekels,
A deep quench approach to the optimal control of
an Allen--Cahn equation with dynamic boundary conditions and double obstacles,
{\it Appl. Math. Optim.} {\bf 71} (2015), 1-24.}
\bibitem{CF}
P. Colli and T. Fukao,
The Allen-Cahn equation with dynamic boundary conditions and mass constraints,
\pier{{\it Math Methods Appl. Sci.}} {\bf 38} (2015), 3950-3967.

\pier{\bibitem{CF2}
P. Colli and T. Fukao, 
Cahn--Hilliard equation with dynamic boundary conditions 
and mass constraint on the boundary, \emph{J. Math. Anal. Appl.}
{\bf 429} (2015), 1190-1213.}

\pier{\bibitem{CGM}
P. Colli, G. Gilardi and G. Marinoschi, 
A boundary control problem for a possibly singular phase 
field system with dynamic boundary conditions, 
{\it J. Math. Anal. Appl.} {\bf 434} (2016), 432-463.}

\bibitem{CGNS_AC}
P. Colli, G. Gilardi, R. Nakayashiki and K. Shirakawa,
A class of quasi-linear Allen-Cahn type equations with dynamic boundary conditions,
{\it Nonlinear Anal.} {\bf 158} (2017), 32-59.

\bibitem{CGSGuid}
P. Colli, G. Gilardi, P. Podio-Guidugli and J. Sprekels, Distributed optimal control
of a nonstandard system of phase field equations, 
{\it \pier{Contin.} Mech. Thermodyn.} {\bf 24} (2012), 437-459.

\bibitem{CGS_Lim}
P. Colli, G. Gilardi and J. Sprekels,
Limiting problems for a nonstandard viscous Cahn-Hilliard 
system with dynamic boundary conditions, \pier{in
{\sl Trends on Applications of Mathematics to Mechanics},}
E. Rocca, U. Stefanelli, L. Truskinovski, A. Visintin (ed.),
{\it Springer INdAM Series} {\bf 27}, Springer, Milan, (2018), 217-242.

\bibitem{CGS_vel}
P. Colli, G. Gilardi and J. Sprekels,
Optimal velocity control of a viscous Cahn-Hilliard system with convection and dynamic boundary conditions, 
{\it SIAM J. Control Optim.\/} {\bf 56} (2018), 1665-1691.

\bibitem{CGS_nonst}
P. Colli, G. Gilardi and J. Sprekels,
Optimal boundary control of a nonstandard viscous \CH~system with dynamic boundary condition,
{\it Nonlinear Anal.\/} {\bf 170} (2018), 171-196.

\bibitem{CGS_Rec}
P. Colli, G. Gilardi and J. Sprekels,
Recent results on the Cahn-Hilliard equation with dynamic boundary conditions,
{\it Vestn. Yuzhno-Ural. Gos. Univ., Ser. Mat. Model. Program.} {\bf 10} (2017), 5-21.

\bibitem{CGS_OPT}
P. Colli, G. Gilardi and J. Sprekels,
A boundary control problem for the viscous \CH~equation 
with dynamic boundary conditions,
{\it Appl. Math. \pier{Optim.}} {\bf 73} (2016), 195-225.

\bibitem{CGSpure}
P. Colli, G. Gilardi and J. Sprekels,
A boundary control problem for the pure \CH~equation
with dynamic boundary conditions,
{\it Adv. Nonlinear Anal.} {\bf 4} (2015), 311-325.

\bibitem{CGS}
P. Colli, G. Gilardi and J. Sprekels,
On the \CH~equation with dynamic 
boundary conditions and a dominating boundary potential,
{\it J. Math. Anal. Appl.\/} {\bf 419} (2014), 972-994.

\bibitem{CGS_nonpf}
P. Colli, G. Gilardi and J. Sprekels, Analysis and boundary 
control of a nonstandard system of phase field equations, 
{\it Milan J. Math.} {\bf 80} (2012), 119-149.

\bibitem{CS_opt}
P. Colli and J. Sprekels,
Optimal boundary control of a nonstandard Cahn-Hilliard system
with dynamic boundary condition and double obstacle inclusions, \pier{in
{\sl Solvability, Regularity, Optimal Control of Boundary Value Problems for PDEs}},
P. Colli, A. Favini, E. Rocca, G. Schimperna, J. Sprekels (ed.),
{\it Springer INdAM Series} {\bf 22}, Springer, Milan, (2017), 151-182.

\pier{\bibitem{FY}
T. Fukao and N. Yamazaki, 
A boundary control problem for the equation and dynamic
boundary condition of Cahn--Hilliard type. \pier{in
{\sl Solvability, Regularity, Optimal Control of Boundary Value Problems for PDEs}},
P. Colli, A. Favini, E. Rocca, G. Schimperna, J. Sprekels (ed.),
{\it Springer INdAM Series} {\bf 22}, Springer, Milan, (2017), 255-280.}

\pier{%
\bibitem{GK}
H. Garcke and P. Knopf,
Weak solutions of the {C}ahn--{H}illiard system with dynamic boundary conditions: a gradient flow approach, preprint arXiv:1810.09817 [math.AP] (2018), 1-27.%
}

\bibitem{GMS_long}
G. Gilardi, A. Miranville and G. Schimperna, 
Long-time behavior of the Cahn-Hilliard equation with irregular potentials and dynamic boundary conditions, {\it Chin. Ann. Math. Ser. B\/}  {\bf 31} (2010), 679-712.

\bibitem{GiMiSchi} 
G. Gilardi, A. Miranville and G. Schimperna,
On the \CH~equation with irregular potentials and dynamic boundary conditions,
{\it Commun. Pure Appl. Anal.\/} {\bf 8} (2009), 881-912.

\bibitem{GS}
G. Gilardi and J. Sprekels,
Asymptotic limits and optimal control for the Cahn-Hilliard 
system with convection and dynamic boundary conditions, 
{\it Nonlinear Analysis} {\bf 178}, (2019), 1-31.

\bibitem{Hint_due}
M. Hinterm\"uller and D. Wegner, Optimal control of a semi-discrete 
Cahn-Hilliard-Navier-Stokes system, 
{\it SIAM J. Control Optim.} {\bf 52} (2014), 747-772.

\bibitem{Hint_uno}
M. Hinterm\"uller and D. Wegner, Distributed optimal control of the Cahn-Hilliard
represent system including the case of a double-obstacle homogeneous free energy density, 
{\it SIAM J. Control Optim.} {\bf 50} (2012), 388-418.

\bibitem{Is} 
H.\ Israel,
Long time behavior of an {A}llen-{C}ahn type equation with a
singular potential and dynamic boundary conditions,
{\it J. Appl. Anal. Comput.} {\bf 2} (2012), 29-56.%

\bibitem{LionsQM}
J.-L. Lions,
Quelques m\'ethodes de r\'esolution des probl\`emes aux limites non lin\'eaires, 
Dunod, Gauthier-Villars, Paris, 1969.

\bibitem{Lions_OPT}
J.-L. Lions,
Contr\^ole optimal de syst\`emes gouverne\'s par des equations aux d\'eriv\'ees partielles,
Dunod, Paris, 1968.


\bibitem{LioMag}
J.-L. Lions and E. Magenes,
Non-homogeneous boundary value problems and applications,
Vol.~I,
Springer, Berlin, 1972.

\andrea{%
\bibitem{LW}
C. Liu and H. Wu,
\pier{An energetic variational approach for the {C}ahn--{H}illiard equation with dynamic boundary conditions: model derivation and mathematical analysis},
{\it Arch. \pier{Ration.} Mech. Anal.} {\bf 233} (2019), 167-247.%
}

\bibitem{Mir_CH}
A. Miranville,
The Cahn-Hilliard equation and some of its variants, 
{\it AIMS Mathematics,} {\bf 2} (2017), 479-544.

\bibitem{MRSS_pf}
A. Miranville, E. Rocca, G. Schimperna and A. Segatti,
The Penrose-Fife phase-field model with coupled dynamic boundary conditions,
\pier{{\it Discrete Contin. Dyn. Syst.} {\bf 34}} (2014), 4259-4290.

\bibitem{MZ} 
A. Miranville and S. Zelik,
\pier{The Cahn-Hilliard equation with singular potentials 
and dynamic boundary conditions,
{\it Discrete Contin. Dyn. Syst.} {\bf 28}} (2010), 275-310.

\bibitem{PRZ} 
J. Pr\"uss, R. Racke and S. Zheng, 
Maximal regularity and asymptotic behavior of solutions 
for the \CH~equation with dynamic boundary conditions,  
{\it Ann. Mat. Pura Appl.~(4)\/}
{\bf 185} (2006), 627-648.

\bibitem{RZ} 
R. Racke and S. Zheng, 
The \CH~equation with dynamic boundary conditions, 
{\it Adv. Differential Equations\/} 
{\bf 8} (2003), 83-110.

\bibitem{RS_OPT}
E. Rocca and J. Sprekels, 
\pier{Optimal distributed control of a nonlocal convective \CH\ equation by the velocity in three dimensions,
{\it SIAM J. Control Optim.}  {\bf 53} (2015), 1654-1680.}

\bibitem{Simon}
J. Simon,
{Compact sets in the space $L^p(0,T; B)$},
{\it Ann. Mat. Pura Appl.~(4)\/} 
{\bf 146} (1987), 65-96.

\bibitem{Trol}
F. Tr\"oltzsch,
Optimal Control of Partial Differential Equations. Theory, Methods and Applications,
{\it Grad. Stud. in Math.,} Vol. {\bf 112}, AMS, Providence, RI, 2010.

\bibitem{WZ} H. Wu and S. Zheng,
Convergence to equilibrium for the \CH~equation with dynamic boundary conditions, 
{\it J. Differential Equations\/}
{\bf 204} (2004), 511-531.

\bibitem{ZL_due}
X. P. Zhao and C. C. Liu, Optimal control for the convective Cahn-Hilliard equation
in 2D case, 
{\it Appl. Math. Optim.} {\bf 70} (2014), 61-82.

\bibitem{ZL}
X. P. Zhao and C. C. Liu, Optimal control of the convective Cahn-Hilliard equation,
{\it Appl. Anal.} {\bf 92} (2013), 1028-1045.

\End{thebibliography}

\End{document}

\bye